\documentclass[reqno]{amsart}
\usepackage{hyperref}

\usepackage{amsmath,amssymb,amsthm}
\usepackage{enumerate}
\usepackage{color}
\usepackage{epsfig}
\usepackage{graphics}
\usepackage{graphicx}
\usepackage{subfigure}

\begin{document}
\title[\hfilneg  \hfil ]{Essential Spectra of Quasi-parabolic Composition Operators on Hardy Spaces of Analytic Functions}
\dedicatory{Dedicated to the memory of Ali Y{\i}ld{\i}z (1976-2006)}

\author[u. g\"{u}l]{u\u{g}ur g\"{u}l}

\address{u\u{g}ur g\"{u}l,  \newline
Sabanc{\i} University, Faculty of Engineering and Natural Sciences,
34956, Tuzla, Istanbul, TURKEY}
\email{\href{mailto:gulugur@gmail.com}{gulugur@gmail.com}}

%\href{mailto:myaddress@wikibooks.org}{myaddress@wikibooks.org}
%\email{{mailto:aozbekler@gmail.com}{aozbekler@gmail.com}}

\thanks{Submitted August 16, 2010}


\begin{abstract}
    In this work we study the essential spectra of composition operators
    on Hardy spaces of analytic functions which might be termed as
    ``quasi-parabolic.'' This is the class of composition operators on
    $H^{2}$ with symbols whose conjugate with the Cayley transform on
    the upper half-plane are of the form $\varphi(z)=$ $z+\psi(z)$, where
    $\psi\in$ $H^{\infty}(\mathbb{H})$ and $\Im(\psi(z)) > \epsilon > 0$.
    We especially examine the case where $\psi$ is discontinuous at
    infinity. A new method is devised to show that this type of
    composition operator fall in a C*-algebra of Toeplitz operators and
    Fourier multipliers. This method enables us to provide new examples
    of essentially normal composition operators and to calculate their
    essential spectra.
\end{abstract}

\maketitle
%\numberwithin{equation}{section}
\newtheorem{theorem}{Theorem}
\newtheorem{acknowledgement}[theorem]{Acknowledgement}
\newtheorem{algorithm}[theorem]{Algorithm}
\newtheorem{axiom}[theorem]{Axiom}
\newtheorem{case}[theorem]{Case}
\newtheorem{claim}[theorem]{Claim}
\newtheorem{conclusion}[theorem]{Conclusion}
\newtheorem{condition}[theorem]{Condition}
\newtheorem{conjecture}[theorem]{Conjecture}
\newtheorem{corollary}[theorem]{Corollary}
\newtheorem{criterion}[theorem]{Criterion}
\newtheorem{definition}[theorem]{Definition}
\newtheorem{example}[theorem]{Example}
\newtheorem{exercise}[theorem]{Exercise}
\newtheorem{lemma}[theorem]{Lemma}
\newtheorem{notation}[theorem]{Notation}
\newtheorem{problem}[theorem]{Problem}
\newtheorem{proposition}[theorem]{Proposition}
\newtheorem{remark}[theorem]{Remark}
\newtheorem{solution}[theorem]{Solution}
\newtheorem{summary}[theorem]{Summary}
\newtheorem*{thma}{Theorem A}
\newtheorem*{thmb}{Theorem B}
\newtheorem*{thmc}{Theorem C}
\newtheorem*{thmd}{Theorem D}
\newcommand{\norm}[1]{\left\Vert#1\right\Vert}
\newcommand{\abs}[1]{\left\vert#1\right\vert}
\newcommand{\set}[1]{\left\{#1\right\}}
\newcommand{\Real}{\mathbb R}
\newcommand{\eps}{\varepsilon}
\newcommand{\To}{\longrightarrow}
\newcommand{\BX}{\mathbf{B}(X)}
\newcommand{\A}{\mathcal{A}}

\section{introduction}

  This work is motivated by the results of Cowen (see [8]) on the spectra of composition operators on $H^{2}(\mathbb{D})$ induced
  by parabolic linear fractional non-automorphisms that fix a point $\xi$ on the boundary. These composition operators are precisely
  the essentially normal linear fractional composition operators [3]. These linear fractional transformations for $\xi= 1$ take
  the form
  \[\varphi_{a}(z)=\frac{2iz+a(1-z)}{2i+a(1-z)}\]
  with $\Im(a) > 0$. Their upper half-plane re-incarnations via the Cayley transform $\mathfrak{C}$ (see page 4) are the translations
\begin{equation*}
  \mathfrak{C}^{-1}\circ\varphi_{a}\circ\mathfrak{C}(w)= w+a
\end{equation*}
   acting on the upper half-plane.

  Cowen [8] has proved that
  \[\sigma(C_{\varphi_{a}})=\sigma_{e}(C_{\varphi_{a}})=\{e^{iat}:t\in [0,\infty)\}\cup\{0\}.\]

Bourdon, Levi, Narayan, Shapiro [3] dealt with composition operators
with symbols $\varphi$ such that the upper half-plane re-incarnation
of $\varphi$ satisfies
\begin{equation*}
\mathfrak{C}^{-1}\circ\varphi\circ\mathfrak{C}(z)=pz+\psi(z) ,
\end{equation*}
 where $p>0$, $\Im(\psi(z))>\epsilon>0$ for all $z\in\mathbb{H}$ and
  $\lim_{z\rightarrow\infty}\psi(z)=$ $\psi_{0}\in$ $\mathbb{H}$
  exist. Their results imply that the essential spectrum of such a composition operator with $p=1$ is
\begin{equation*}
  \{e^{i\psi_{0}t}:t\in [0,\infty)\}\cup\{0\}.
\end{equation*}
In this work we are interested in composition operators whose
symbols $\varphi$ have upper half-plane re-incarnation
\begin{equation*}
\mathfrak{C}^{-1}\circ\varphi\circ\mathfrak{C}(z)=z+\psi(z)
\end{equation*}
 for a bounded analytic function
  $\psi$ satisfying $\Im(\psi(z)) > \epsilon > 0$ for all $z\in$
  $\mathbb{H}$. This class will obviously include those studied in [3] with $p=1$. However we will be particularly interested in the case where $\psi$ does not have a limit at infinity.
  We call such composition operators ``quasi-parabolic.'' Our
  most precise result is obtained when the boundary values of $\psi$
  lie in $QC$, the space of quasi-continuous functions on
  $\mathbb{T}$, which is defined as
  \[QC = [H^{\infty}+C(\mathbb{T})]\cap[\overline{H^{\infty}+C(\mathbb{T})}].\]

  We recall that the set of cluster points $\mathcal{C}_{\xi}(\psi)$ of $\psi\in
  H^{\infty}$ is defined to be the set of points
  $z\in$ $\mathbb{C}$ for which there is a sequence $\{z_{n}\}\subset$
  $\mathbb{D}$ so that $z_{n}\rightarrow$ $\xi$ and
  $\psi(z_{n})\rightarrow$ $z$.

  In particular we prove the following theorem.

\begin{thmb}
Let $\varphi:\mathbb{D}\rightarrow$ $\mathbb{D}$ be an analytic
self-map of $\mathbb{D}$ such that
$$\varphi(z)=\frac{2iz+\eta(z)(1-z)}{2i+\eta(z)(1-z)} , $$
where $\eta\in$ $H^{\infty}(\mathbb{D})$ with $\Im(\eta(z)) >
\epsilon> 0$ for all $z\in$ $\mathbb{D}$. If $\eta\in$ $QC\cap
H^{\infty}$ then we have
\begin{itemize}
\item\quad $C_{\varphi}:$ $H^{2}(\mathbb{D})\rightarrow$ $H^{2}(\mathbb{D})$ is essentially normal
\item\quad$\sigma_{e}(C_{\varphi})=\{e^{izt}:t\in [0,\infty],z\in\mathcal{C}_{1}(\eta)\}\cup\{0\}$,
\end{itemize}
 where\quad $\mathcal{C}_{1}(\eta)$ is the set of cluster points of $\eta$ at $1$.
\end{thmb}

Moreover, for general $\eta\in H^{\infty}$ with
$\Im(\eta(z))>\epsilon> 0$ (but no requirement that $\eta\in QC$),
we have
\begin{equation*}
\sigma_{e}(C_{\varphi})\supseteq\{e^{izt}:t\in [0,\infty),
    z\in\mathcal{R}_{1}(\eta)\}\cup\{0\} ,
\end{equation*}
where the local essential
  range $\mathcal{R}_{\xi}(\eta)$ of an $\eta\in$
  $L^{\infty}(\mathbb{T})$ at $\xi\in$ $\mathbb{T}$ is defined
  to be the set of points $z\in$ $\mathbb{C}$ so that, for all $\varepsilon > 0$ and $\delta > 0$, the set
\begin{equation*}
  \eta^{-1}(B(z,\varepsilon))\cap\{e^{it}:\mid t-t_{0}\mid\leq\delta\}
\end{equation*}
  has positive Lebesgue measure, where $e^{it_{0}}=$ $\xi$. We
  note that ([30]) for functions $\eta\in$ $QC\cap H^{\infty}$,
\begin{equation*}
  \mathcal{R}_{\xi}(\eta)=\mathcal{C}_{\xi}(\eta).
\end{equation*}
   The local essential range $\mathcal{R}_{\infty}(\psi)$ of $\psi\in
L^{\infty}(\mathbb{R})$ at $\infty$
  is defined as the set of points $z\in$ $\mathbb{C}$ so that, for all
$\varepsilon > 0$ and $n > 0$, we have
\begin{equation*}
    \lambda(\psi^{-1}(B(z,\varepsilon))\cap(\mathbb{R}-[-n,n]))>0 ,
\end{equation*}
where $\lambda$ is the Lebesgue measure on $\mathbb{R}$.

   The Cayley transform induces a natural isometric isomorphism between $H^{2}(\mathbb{D})$ and $H^{2}(\mathbb{H})$. Under this identification ``quasi-parabolic" composition operators correspond to operators of the form
\begin{equation*}
  T_{\frac{\varphi(z)+i}{z+i}}C_{\varphi}=C_{\varphi}+T_{\frac{\psi(z)}{z+i}}C_{\varphi} ,
\end{equation*}
   where $\varphi(z)=$ $z+\psi(z)$ with $\psi\in$ $H^{\infty}$ on the upper half-plane, and $T_{\phi}$ is the multiplication operator by
   $\phi$.

   We work on the upper half-plane and use Banach algebra
   techniques to compute the essential spectra of operators that correspond to ``quasi-parabolic"
   operators. Our treatment is motivated by [12] where the
   translation operators are considered as Fourier multipliers on $H^{2}$ (we refer the reader to [22] for the definition and properties of Fourier multipliers).
   Throughout the present work, $H^{2}(\mathbb{H})$ will be considered as a
   closed subspace of $L^{2}(\mathbb{R})$ via the boundary values.
   With the help of Cauchy integral formula we prove an integral
   formula that gives composition operators as integral operators.
   Using this integral formula we show that operators that correspond to``quasi-parabolic"
   operators fall in a C*-algebra generated by Toeplitz operators
   and Fourier multipliers.

The remainder of this paper is organized as follows. In section 2 we
give the basic definitions and preliminary material that we will use
throughout. For the benefit of the reader we explicitly recall some
facts from Banach algebras and Toeplitz operators. In Section 3 we
first prove an integral representation formula for composition
operators on $H^{2}$ of the upper half-plane. Then we use this
integral formula to prove that a ``quasi-parabolic" composition
operator is written as a series of Toeplitz operators and Fourier
multipliers which converges in operator norm. In section 4 we
analyze the C*-algebra generated by Toeplitz operators with
$QC(\mathbb{R})$ symbols and Fourier multipliers modulo compact
operators. We show that this C*-algebra is commutative and we
identify its maximal ideal space using a related theorem of Power
(see [23]). In section 5, using the machinery developed in sections
3 and 4, we determine the essential spectra of "quasi-parabolic"
composition operators. We also give an example of a
``quasi-parabolic" composition operator $C_{\varphi}$ for which
$\psi\in$ $QC(\mathbb{R})$ but does not have a limit at infinity and
compute its essential spectrum.

In the last section we examine the case of $C_{\varphi}$ with
\begin{equation*}
    \varphi(z)=z+\psi(z) ,
\end{equation*}
where $\psi\in H^{\infty}(\mathbb{H})$,\quad $\Im(\psi(z))>\epsilon>0$
but $\psi$ is not necessarily in $QC(\mathbb{R})$. Using Power's theorem on the C*-algebra generated by Toeplitz operators with $L^{\infty}(\mathbb{R})$ symbols and Fourier multipliers, we prove the result
\begin{equation*}
    \sigma_{e}(C_{\varphi})\supseteq\{e^{izt}:z\in\mathcal{R}_{\infty}(\psi),t\in [0,\infty)\}\cup\{0\} ,
\end{equation*}
where $\varphi(z)=z+\psi(z)$,\quad $\psi\in$ $H^{\infty}$ with $\Im(\psi(z)) > \epsilon > 0$.

\section{Notation and Preliminaries}

In this section we fix the notation that we will use throughout and
recall some preliminary facts that will be used in the sequel.

Let $S$ be a compact Hausdorff topological space. The space of all
complex valued continuous functions on $S$ will be denoted by
$C(S)$. For any $f\in C(S)$, $\parallel f\parallel_{\infty}$ will
denote the sup-norm of $f$, i.e. $$\parallel
f\parallel_{\infty}=\sup\{\mid f(s)\mid:s\in S\}.$$ For a Banach
space $X$, $K(X)$ will denote the space of all compact operators on
$X$ and $\mathcal{B}(X)$  will denote the space of all bounded
linear operators on $X$. The open unit disc will be denoted by
$\mathbb{D}$, the open upper half-plane will be denoted by
$\mathbb{H}$, the real line will be denoted by $\mathbb{R}$ and the
complex plane will be denoted by $\mathbb{C}$. The one point
compactification of $\mathbb{R}$ will be denoted by
$\dot{\mathbb{R}}$ which is homeomorphic to $\mathbb{T}$. For any
$z\in$ $\mathbb{C}$, $\Re(z)$ will denote the real part, and
$\Im(z)$ will denote the imaginary part of $z$, respectively. For
any subset $S\subset$ $B(H)$, where $H$ is a Hilbert space, the
C*-algebra generated by $S$ will be denoted by $C^{*}(S)$. The
Cayley transform $\mathfrak{C}$ will be defined by
\begin{equation*}
\mathfrak{C}(z)=\frac{z-i}{z+i}.
\end{equation*}
 For any $a\in$ $L^{\infty}(\mathbb{R})$ (or $a\in$
$L^{\infty}(\mathbb{T})$), $M_{a}$ will be the multiplication
operator on $L^{2}(\mathbb{R})$ (or $L^{2}(\mathbb{T})$) defined as
\begin{equation*}
M_{a}(f)(x)=a(x)f(x).
\end{equation*}
For convenience, we remind the reader of the rudiments of Gelfand
theory of commutative Banach algebras and Toeplitz operators.

Let $A$ be a commutative Banach algebra. Then its maximal ideal space
$M(A)$ is defined as
\begin{equation*}
    M(A)=\{x\in A^{*}:x(ab)=x(a)x(b)\quad\forall a,b\in A\}
\end{equation*}
where $A^{*}$ is the dual space of $A$. If $A$ has identity then
$M(A)$ is a compact Hausdorff topological space with the weak*
topology. The Gelfand transform $\Gamma:A\rightarrow C(M(A))$ is
defined as
\begin{equation*}
    \Gamma(a)(x)=x(a).
\end{equation*}
 If $A$ is a commutative C*-algebra with
identity, then $\Gamma$ is an isometric *-isomorphism between $A$
and $C(M(A))$. If $A$ is a C*-algebra and $I$ is a two-sided closed
ideal of $A$, then the quotient algebra $A/I$ is also a C*-algebra
(see [1] and [10]).
 For $a\in A$ the spectrum $\sigma_{A}(a)$ of $a$ on $A$
is defined as
\begin{equation*}
    \sigma_{A}(a)=\{\lambda\in\mathbb{C}:\lambda e-a\ \ \textrm{is not invertible in}\ A\},
\end{equation*}
where $e$ is the identity of $A$. We will
use the spectral permanency property of C*-algebras (see [25], pp.
283 and [10], pp.15); i.e. if $A$ is a C*-algebra with identity and
$B$ is a closed *-subalgebra of $A$, then for any $b\in B$ we have
\begin{equation}
\sigma_{B}(b)=\sigma_{A}(b).
\end{equation}
To compute essential spectra we employ the following important fact
(see [25], pp. 268 and [10], pp. 6, 7): If $A$ is a commutative Banach
algebra with identity then for any $a\in A$ we have
\begin{equation}
    \sigma_{A}(a)=\{\Gamma(a)(x)=x(a):x\in M(A)\}.
\end{equation}
In general (for $A$ not necessarily commutative), we have
\begin{equation}
    \sigma_{A}(a)\supseteq\{x(a):x\in M(A)\}.
\end{equation}

For a Banach algebra $A$, we denote by $com(A)$ the closed ideal in
$A$ generated by the commutators
$\{a_{1}a_{2}-a_{2}a_{1}:a_{1},a_{2}\in A\}$. It is an algebraic
fact that the quotient algebra $A/com(A)$ is a commutative Banach
algebra. The reader can find detailed information about Banach and
C*-algebras in [25] and [10] related to what we have reviewed so
far.

The essential spectrum $\sigma_{e}(T)$ of an operator $T$ acting on
a Banach
  space $X$ is the spectrum of the coset of $T$ in the Calkin algebra
  $\mathcal{B}(X)/K(X)$, the algebra of bounded linear operators modulo
  compact operators. The well known Atkinson's theorem identifies the essential
  spectrum of $T$ as the set of all $\lambda\in$ $\mathbb{C}$ for
  which $\lambda I-T$ is not a Fredholm operator. The essential norm of $T$ will be denoted by $\parallel T\parallel_{e}$ which is defined as
\begin{equation*}
 \parallel T\parallel_{e}=\inf\{\parallel T+K\parallel:K\in K(X)\}
\end{equation*}
   The bracket $[\cdot]$ will denote the equivalence class modulo
  $K(X)$. An operator $T\in\mathcal{B}(H)$ is called essentially
  normal if $T^{*}T-TT^{*}\in K(H)$ where $H$ is a Hilbert space and
  $T^{*}$ denotes the Hilbert space adjoint of $T$.

  For $1\leq p < \infty$ the Hardy space of the unit disc will be
denoted by $H^{p}(\mathbb{D})$ and the Hardy space of the upper
half-plane will be denoted by $H^{p}(\mathbb{H})$.

  The two Hardy spaces $H^{2}(\mathbb{D})$ and $H^{2}(\mathbb{H})$
    are isometrically isomorphic. An isometric isomorphism $\Phi:H^{2}(\mathbb{D})\longrightarrow$ $H^{2}(\mathbb{H})$
is given by
   \begin{equation}\label{nice1}
    \Phi(g)(z)=
   \bigg(\frac{1}{\sqrt{\pi}(z+i)}\bigg)g\bigg(\frac{z-i}{z+i}\bigg)
   \end{equation}
The mapping $\Phi$ has an inverse
$\Phi^{-1}:H^{2}(\mathbb{H})\longrightarrow$ $H^{2}(\mathbb{D})$
given by
\begin{equation*}\label{nice2}
\Phi^{-1}(f)(z)= \frac{e^\frac{i\pi}{2}(4\pi)^\frac{1}{2}}{(1-z)}
f\bigg(\frac{i(1+z)}{1-z}\bigg)
\end{equation*}
For more details see [15, pp. 128-131] and [19].

Using the isometric isomorphism $\Phi$, one may transfer Fatou's
theorem in the unit disc case to upper half-plane and may embed
   $H^{2}(\mathbb{H})$ in $L^{2}(\mathbb{R})$ via
   $f\longrightarrow$ $f^{*}$ where $f^{*}(x)=$
   $\lim_{y\rightarrow 0}f(x+iy)$. This embedding is an
   isometry.

Throughout the paper, using $\Phi$, we will go back and forth
between $H^{2}(\mathbb{D})$ and $H^{2}(\mathbb{H})$. We use the
property that $\Phi$ preserves spectra, compactness and essential
spectra i.e. if $T\in\mathcal{B}(H^{2}(\mathbb{D}))$ then
\begin{equation*}
    \sigma_{\mathcal{B}(H^{2}(\mathbb{D}))}(T)=\sigma_{\mathcal{B}(H^{2}(\mathbb{H}))}(\Phi\circ
    T\circ\Phi^{-1}),
\end{equation*}
$K\in K(H^{2}(\mathbb{D}))$ if and only if $\Phi\circ
K\circ\Phi^{-1} \in K(H^{2}(\mathbb{H}))$ and hence we have
\begin{equation}
    \sigma_{e}(T)=\sigma_{e}(\Phi\circ T\circ\Phi^{-1}).
\end{equation}
We also note that $T\in\mathcal{B}(H^{2}(\mathbb{D}))$ is
essentially normal if and only if $\Phi\circ
T\circ\Phi^{-1}\in\mathcal{B}(H^{2}(\mathbb{H}))$ is essentially
normal.

The Toeplitz operator with symbol $a$ is defined as
    $$T_{a}=P M_{a}|_{H^{2}} ,$$
where $P$ denotes the orthogonal
projection of $L^{2}$ onto $H^{2}$. A good reference about Toeplitz
operators on $H^{2}$ is Douglas' treatise ([11]). Although the
Toeplitz operators treated in [11] act on the
 Hardy space of the unit disc, the results can be transfered
 to the upper half-plane case using the isometric isomorphism $\Phi$
 introduced by equation (4). In the sequel the following identity
 will be used:
\begin{equation}
     \Phi^{-1}\circ T_{a}\circ\Phi=T_{a\circ \mathfrak{C}^{-1}} ,
\end{equation}
where $a\in L^{\infty}(\mathbb{R})$. We also employ the fact
\begin{equation}
    \parallel T_{a}\parallel_{e}=\parallel
    T_{a}\parallel=\parallel a\parallel_{\infty}
\end{equation}
 for any $a\in L^{\infty}(\mathbb{R})$, which is a consequence
 of Theorem 7.11 of [11] (pp. 160--161) and equation (6). For any subalgebra $A\subseteq L^{\infty}(\mathbb{R})$ the Toeplitz C*-algebra generated
 by symbols in $A$ is defined to be
\begin{equation*}
 \mathcal{T}(A)=C^{*}(\{T_{a}:a\in A\}).
\end{equation*}
 It is a well-known result of Sarason (see [26],[28] and also [24]) that the
 set of functions
\begin{equation*}
 H^{\infty}+C=\{f_{1}+f_{2}:f_{1}\in H^{\infty}(\mathbb{D}),f_{2}\in C(\mathbb{T})\}
\end{equation*}
 is a closed subalgebra of $L^{\infty}(\mathbb{T})$. The following theorem of
 Douglas [11] will be used in the sequel.
\begin{theorem} [\scshape Douglas' Theorem]
    Let $a$,$b\in$ $H^{\infty}+C$ then the semi-commutators
    $$T_{ab}-T_{a}T_{b}\in K(H^{2}(\mathbb{D})),\quad T_{ab}-T_{b}T_{a}\in K(H^{2}(\mathbb{D})), $$
    and hence the commutator
    $$[T_{a},T_{b}]=T_{a}T_{b}-T_{b}T_{a}\in K(H^{2}(\mathbb{D})).$$ \label{thmDouglas}
\end{theorem}
Let $QC$ be the C*-algebra of functions in $H^{\infty}+C$ whose
complex conjugates also belong to $H^{\infty}+C$. Let us also define
the upper half-plane version of $QC$ as the following:
\begin{equation*}
    QC(\mathbb{R})=\{a\in L^{\infty}(\mathbb{R}):a\circ\mathfrak{C}^{-1}\in QC\}.
\end{equation*}
Going back and forth with Cayley transform one can deduce that
$QC(\mathbb{R})$ is a closed subalgebra of $L^{\infty}(\mathbb{R})$.

By Douglas' theorem and equation (6), if $a$, $b\in QC(\mathbb{R})$,
then
\begin{equation*}
T_{a}T_{b}-T_{ab}\in K(H^{2}(\mathbb{H})).
\end{equation*}
Let $scom(QC(\mathbb{R}))$ be the closed ideal in
$\mathcal{T}(QC(\mathbb{R}))$ generated by the semi-commutators
$\{T_{a}T_{b}-T_{ab}:a, b\in QC(\mathbb{R})\}$. Then we have
\begin{equation*}
    com(\mathcal{T}(QC(\mathbb{R})))\subseteq
    scom(QC(\mathbb{R}))\subseteq K(H^{2}(\mathbb{H})) .
\end{equation*}
 By Proposition 7.12 of [11] and equation (6) we have
\begin{equation}
    com(\mathcal{T}(QC(\mathbb{R})))=scom(QC(\mathbb{R}))=K(H^{2}(\mathbb{H})) .
\end{equation}
 Now consider the symbol map
\begin{equation*}
 \Sigma:QC(\mathbb{R})\rightarrow\mathcal{T}(QC(\mathbb{R}))
\end{equation*}
  defined as
 $\Sigma(a)=T_{a}$. This map is linear but not necessarily multiplicative; however if we let $q$ be
 the quotient map
\begin{equation*}
    q:\mathcal{T}(QC(\mathbb{R})) \rightarrow \mathcal{T}(QC(\mathbb{R}))/scom(QC(\mathbb{R})) ,
\end{equation*}
then $q\circ\Sigma$ is multiplicative; moreover by equations (7) and
(8), we conclude that
 $q\circ\Sigma$ is an isometric *-isomorphism from $QC(\mathbb{R})$
 onto $\mathcal{T}(QC(\mathbb{R}))/K(H^{2}(\mathbb{H}))$.

\begin{definition}
Let
$\varphi:\mathbb{D}\longrightarrow$ $\mathbb{D}$ or
$\varphi:\mathbb{H}\longrightarrow$ $\mathbb{H}$ be a holomorphic
self-map of the unit disc or the upper half-plane. The
\emph{composition
 operator} $C_{\varphi}$ on $H^{p}(\mathbb{D})$ or $H^{p}(\mathbb{H})$ with symbol $\varphi$ is defined by
\begin{equation*}
C_{\varphi}(g)(z)= g(\varphi(z)),\qquad
z\in\mathbb{D}\quad\textrm{or}\quad z\in\mathbb{H}.
\end{equation*}
\end{definition}

 Composition operators of the unit disc are always
 bounded [9] whereas composition operators of the upper half-plane are not
 always bounded. For the boundedness problem of composition operators of the upper half-plane see
 [19].

The composition operator $C_{\varphi}$ on $H^{2}(\mathbb{D})$ is
carried over to
$(\frac{\tilde{\varphi}(z)+i}{z+i})C_{\tilde{\varphi}}$ on
$H^{2}(\mathbb{H})$ through $\Phi$, where $\tilde{\varphi}=$
$\mathfrak{C}\circ\varphi\circ\mathfrak{C}^{-1}$, i.e. we have
\begin{equation}
    \Phi C_{\varphi}\Phi^{-1} =
    T_{(\frac{\tilde{\varphi}(z)+i}{z+i})}C_{\tilde{\varphi}}.
\end{equation}

However this gives us the boundedness of
 $C_{\varphi}:H^{2}(\mathbb{H})$ $\rightarrow$ $H^{2}(\mathbb{H})$ for
\begin{equation*}
    \varphi(z)=pz+\psi(z) ,
\end{equation*}
 where $p > 0$, $\psi\in$ $H^{\infty}$ and $\Im(\psi(z))>\epsilon>0$ for all $z\in\mathbb{H}$:

Let $\tilde{\varphi}:\mathbb{D}\rightarrow$ $\mathbb{D}$ be an analytic
self-map of $\mathbb{D}$ such that $\varphi=$ $\mathfrak{C}^{-1}\circ\tilde{\varphi}\circ\mathfrak{C}$,
then we have
\begin{equation*}
    \Phi C_{\tilde{\varphi}}\Phi^{-1}=T_{\tau}C_{\varphi}
\end{equation*}
where
\begin{equation*}
    \tau(z)=\frac{\varphi(z)+i}{z+i}.
\end{equation*}
If
\begin{equation*}
    \varphi(z)=pz+\psi(z)
\end{equation*}
with $p > 0$, $\psi\in H^{\infty}$ and $\Im(\psi(z)) > \epsilon> 0$, then
$T_{\frac{1}{\tau}}$ is a bounded operator. Since $\Phi C_{\tilde{\varphi}}\Phi^{-1}$  is always bounded we conclude
that $C_{\varphi}$ is bounded on $H^{2}(\mathbb{H})$.

We recall that any function in $H^{2}(\mathbb{H})$ can be recovered
from its boundary values by means of the Cauchy integral. In fact we
have [16, pp. 112--116] if $f\in$ $H^{2}(\mathbb{H})$ and if $f^{*}$
is its non-tangential boundary value function on $\mathbb{R}$, then
\begin{equation}
    f(z)= \frac{1}{2\pi i}\int_{-\infty}^{+\infty}\frac{f^{*}(x)dx}{x-z},\qquad
    z\in\mathbb{H}.
\end{equation}

The Fourier transform $\mathcal{F}f$ of $f\in$
$\mathcal{S}(\mathbb{R})$ (the Schwartz space, for a definition see
[25, sec. 7.3, pp.~168] and [33, pp.~134]) is defined by
\begin{equation*}
(\mathcal{F}f)(t)=\frac{1}{\sqrt{2\pi}}\int_{-\infty}^{+\infty}e^{-itx}f(x)dx.
\end{equation*}
The Fourier transform extends to an invertible isometry from
$L^{2}(\mathbb{R})$ onto itself with inverse
\begin{equation*}
(\mathcal{F}^{-1}f)(t)=\frac{1}{\sqrt{2\pi}}\int_{-\infty}^{+\infty}
e^{itx}f(x)dx.
\end{equation*}
The following is a consequence of a theorem due to Paley and Wiener
[16, pp.~110--111]. Let $1 < p < \infty$. For $f\in$
$L^{p}(\mathbb{R})$, the following assertions are equivalent:
\begin{enumerate}
    \item[($i$)]  $f\in$ $H^{p}$,
    \item[($ii$)] $\textrm{supp}(\hat{f})\subseteq$ $[0,\infty)$
\end{enumerate}

A reformulation of the Paley-Wiener theorem says that the image of
$H^{2}(\mathbb{H})$ under the Fourier transform is
$L^{2}([0,\infty))$.

 By the Paley-Wiener theorem we observe that the operator
$$D_{\vartheta}=\mathcal{F}^{-1}M_{\vartheta}\mathcal{F}$$
for $\vartheta\in C([0,\infty])$ maps $H^{2}(\mathbb{H})$ into
itself, where $C([0,\infty])$ denotes the set of continuous
functions on $[0,\infty)$ which have limits at infinity. Since
$\mathcal{F}$ is unitary we also observe that
\begin{equation}
    \parallel D_{\vartheta}\parallel=\parallel M_{\vartheta}\parallel=\parallel\vartheta\parallel_{\infty}
\end{equation}
Let $F$ be defined as
\begin{equation}
F =\{D_{\vartheta}\in B(H^{2}(\mathbb{H})):\vartheta\in
C([0,\infty])\} .
\end{equation}
We observe that $F$ is a commutative C*-algebra with identity and
the map $D:C([0,\infty])\rightarrow F$ given by
\begin{equation*}
D(\vartheta)=D_{\vartheta}
\end{equation*}
is an isometric *-isomorphism by equation (11). Hence $F$ is
isometrically *-isomorphic to $C([0,\infty])$. The operator
$D_{\vartheta}$ is usually called a ``Fourier Multiplier.''

An important example of a Fourier multiplier is the translation
operator $S_{w}:H^{2}(\mathbb{H})\rightarrow H^{2}(\mathbb{H})$
defined as
\begin{equation*}
S_{w}f(z)=f(z+w)
\end{equation*}
where $w\in\mathbb{H}$. We recall that
\begin{equation*}
S_{w}=D_{\vartheta}
\end{equation*}
where $\vartheta(t)=e^{iwt}$ (see [12] and [14]). Other examples of
Fourier multipliers that we will need come from convolution
operators defined in the following way:
\begin{equation}
    K_{n}f(x)=\frac{1}{2\pi
    i}\int_{-\infty}^{\infty}\frac{-f(w)dw}{(x-w+i\alpha)^{n+1}} ,
\end{equation}
where $\alpha\in\mathbb{R}^{+}$. We observe that
\begin{eqnarray*}
    & &\mathcal{F}K_{n}f(x)=\int_{-\infty}^{\infty}e^{-itx}\bigg(\int_{-\infty}^{\infty}\frac{-f(w)dw}{(t-w+i\alpha)^{n+1}}\bigg)dt\\
    =& &\int_{-\infty}^{\infty}\int_{-\infty}^{\infty}\frac{e^{-i(t-w)}e^{-iwx}(-f(w))}{(t-w+i\alpha)^{n+1}}dwdt\\
    =& &\bigg(\int_{-\infty}^{\infty}\frac{-e^{-ivx}dv}{(v+i\alpha)^{n+1}}\bigg)\bigg(\int_{-\infty}^{\infty}e^{-iwx}f(w)dw\bigg).
\end{eqnarray*}
 Since
\begin{equation*}
    \int_{-\infty}^{\infty}\frac{-e^{-ivx}dv}{(v+i\alpha)^{n+1}}=\frac{(-ix)^{n}e^{-\alpha x}}{n!} ,
\end{equation*}
this implies that
\begin{equation}
    K_{n}=D_{\vartheta_{n}}
\end{equation}
where
\begin{equation*}
\vartheta_{n}(t)=\frac{(-it)^{n}e^{-\alpha t}}{n!}.
\end{equation*}
 For $p>0$ the dilation operator $V_{p}\in\mathcal{B}(H^{2}(\mathbb{H}))$
is defined as
\begin{equation}
    V_{p}f(z)=f(pz).
\end{equation}

\section{an approximation scheme for composition operators on hardy spaces of the upper half-plane}

In this section we devise an integral representation formula for composition operators and using this integral formula we develop an approximation scheme for composition
operators induced by maps of the form
\begin{equation*}
    \varphi(z)=pz+\psi(z) ,
\end{equation*}
where $p > 0$ and $\psi\in$ $H^{\infty}$ such that $\Im(\psi(z))>\epsilon>0$ for all $z\in\mathbb{H}$. By the preceding section we know that these
maps induce bounded composition operators on $H^{2}(\mathbb{H})$.
We approximate these operators by linear combinations of Toeplitz operators
and Fourier multipliers.  In establishing
this approximation scheme our main tool is the integral representation formula that
we prove below.

One can use equation (10) to represent composition operators with an
integral kernel under some conditions on the analytic symbol
$\varphi:\mathbb{H}\rightarrow\mathbb{H}$. One may apply the
argument (using the Cayley transform) done after equation (4) to
$H^{\infty}(\mathbb{H})$ to show that
\begin{equation*}
    \lim_{t\rightarrow 0}\varphi(x+it)=\varphi^{*}(x)
\end{equation*}
     exists for
almost every $x\in$ $\mathbb{R}$. The most important condition that
we will impose on $\varphi$ is $\Im(\varphi^{*}(x)) > 0$ for almost
every $x\in$ $\mathbb{R}$. We have the following proposition.

\begin{proposition} Let $\varphi:\mathbb{H}\rightarrow$
    $\mathbb{H}$ be an analytic function such that the non-tangential
    boundary value function $\varphi^{*}$ of $\varphi$ satisfies
    $\Im(\varphi^{*}(x)) > 0$ for almost every $x\in$ $\mathbb{R}$.
    Then the composition operator $C_{\varphi}$ on
    $H^{2}(\mathbb{H})$ is given by
\begin{equation*}
(C_{\varphi}f)^{*}(x)=\frac{1}{2\pi
i}\int_{-\infty}^{\infty}\frac{f^{*}(\xi)d\xi}{\xi-\varphi^{*}(x)}\quad
\text{for almost every}\quad x\in\mathbb{R}.
\end{equation*}
\end{proposition}

\begin{proof}
    By equation (10) above one has
    \begin{equation*}\label{nice6}
        C_{\varphi}(f)(x+it)= \frac{1}{2\pi
        i}\int_{-\infty}^{\infty}\frac{f^{*}(\xi)d\xi}{\xi-\varphi(x+it)}.
    \end{equation*}
    Let $x\in$ $\mathbb{R}$ be such that $\lim_{t\rightarrow
    0}\varphi(x+it)=$ $\varphi^{*}(x)$ exists and $\Im(\varphi^{*}(x))> 0$. We have
    \begin{eqnarray}
        &\nonumber &\left\lvert C_{\varphi}(f)(x+it)-\frac{1}{2\pi i}\int_{-\infty}^{\infty}\frac{f^{*}(\xi)d\xi}{\xi-\varphi^{*}(x)}\right\rvert\\
        =& &\left\lvert\frac{1}{2\pi i}\int_{-\infty}^{\infty}\frac{f^{*}(\xi)d\xi}{\xi-\varphi(x+it)}-\frac{1}{2\pi i}\int_{-\infty}^{\infty}\frac{f^{*}(\xi)d\xi}{\xi-\varphi^{*}(x)}\right\rvert \\
        =&\nonumber &\frac{1}{2\pi}\mid\varphi(x+it)-\varphi^{*}(x)\mid \left\lvert\int_{-\infty}^{\infty}\frac{f^*(\xi)d\xi}{(\xi-\varphi(x+it))(\xi-\varphi^{*}(x))}\right\rvert \\
        \leq &\nonumber &\frac{\mid\varphi(x+it)-\varphi^{*}(x)\mid}{2\pi}\| f\|_{2}\bigg(\int_{-\infty}^{\infty}\frac{d\xi}{(\mid(\xi-\varphi(x+it))(\xi-\varphi^{*}(x))\mid)^2}\bigg)^\frac{1}{2} ,
    \end{eqnarray}
   by Cauchy-Schwarz inequality. When
    $\mid\varphi(x+it)-\varphi^{*}(x)\mid < \varepsilon$, by triangle
    inequality, we have
    \begin{equation}
        \mid\xi-\varphi(x+it)\mid\geq\mid\xi-\varphi^{*}(x)\mid-\varepsilon .
    \end{equation}

    Fix $\varepsilon_{0} > 0$ such that
\begin{equation*}
        \varepsilon_{0}=\frac{\inf\{\mid\xi-\varphi^{*}(x)\mid:\xi\in\mathbb{R}\}}{2}.
\end{equation*}
    This is possible since $\Im(\varphi^{*}(x)) > 0$.

    Choose $\varepsilon  > 0$ such that $\varepsilon_{0} >$
    $\varepsilon$. Since $\lim_{t\rightarrow 0}\varphi(x+it)=$
    $\varphi^{*}(x)$ exists, there exists $\delta  > 0$ such that for all
    $0<t<\delta$ we have
\begin{equation*}
    \mid\varphi(x+it)-\varphi^{*}(x)\mid <\varepsilon <\varepsilon_{0}.
\end{equation*}
     So by equation (17) one has
    \begin{equation}
         \mid\xi-\varphi(x+it)\mid\geq
        \mid\xi-\varphi^{*}(x)\mid-\varepsilon_{0}\geq\varepsilon_{0}
    \end{equation}
    for all $t$ such that $0<t<\delta$. By equation (18) we have
    \begin{equation*}
         \frac{1}{\mid\xi-\varphi(x+it)\mid}\leq
        \frac{1}{\mid\xi-\varphi^{*}(x)\mid-\varepsilon_{0}}.
    \end{equation*}
    which implies that
    \begin{eqnarray}
        &\nonumber &\int_{-\infty}^{\infty}\frac{d\xi}{(\mid(\xi-\varphi(x+it))(\xi-\varphi^{*}(x))\mid)^2}\\
        \leq& &\int_{-\infty}^{\infty}\frac{d\xi}{\mid\xi-\varphi^{*}(x)\mid^{4}-\varepsilon_{0}\mid\xi-\varphi^{*}(x)\mid^{2}}
    \end{eqnarray}
    By the right hand-side inequality of equation (18), the integral on the right-hand side of equation
    (19) converges and its value only depends on $x$ and
    $\varepsilon_{0}$. Let $M_{\varepsilon_{0},x}$ be the value of
    that integral, then by equations (16) and (19) we have
    \begin{eqnarray*}
        & &\left\lvert C_{\varphi}(f)(x+it)-\frac{1}{2\pi i}\int_{-\infty}^{\infty}\frac{f^{*}(\xi)d\xi}{\xi-\varphi^{*}(x)}\right\rvert\\
        \leq& &\frac{\mid\varphi(x+it)-\varphi^{*}(x)\mid}{2\pi}\|f\|_{2}\bigg(\int_{-\infty}^{\infty}\frac{d\xi}{\mid\xi-\varphi^{*}(x)\mid^{4}-\varepsilon_{0}\mid\xi-\varphi^{*}(x)\mid^{2}}\bigg)^{\frac{1}{2}}\\
        =& &\frac{\mid\varphi(x+it)-\varphi^{*}(x)\mid}{2\pi}\|f\|_{2}(M_{\varepsilon_{0},x})^{\frac{1}{2}}\leq\frac{\varepsilon}{2\pi}\|f\|_{2}(M_{\varepsilon_{0},x})^{\frac{1}{2}}.
    \end{eqnarray*}
    Hence we have
\begin{equation*}
    \lim_{t\rightarrow 0}C_{\varphi}(f)(x+it)= C_{\varphi}(f)^{*}(x)= \frac{1}{2\pi i}\int_{-\infty}^{\infty}\frac{f^{*}(\xi)d\xi}{\xi-\varphi^{*}(x)}
\end{equation*}
    for $x\in$ $\mathbb{R}$ almost everywhere.
\end{proof}

Throughout the rest of the paper we will identify a function $f$ in
$H^{2}$ or $H^{\infty}$ with its boundary function $f^{*}$.

 We continue with the following simple geometric lemma that will be helpful in our
 task.

 \begin{lemma}
    Let $K\subset$ $\mathbb{H}$ be a compact
    subset of $\mathbb{H}$. Then there is an $\alpha\in$ $\mathbb{R}^{+}$
    such that $\sup\{\mid\frac{\alpha i-z}{\alpha}\mid:z\in K\}<$
    $\delta < 1$ for some $\delta\in$ $(0,1)$.
    \label{lem1}
\end{lemma}

\begin{proof}
    Let $\varepsilon=$ $\inf\{\Im(z):z\in K\}$,
    $R_{1}=$ $\sup\{\Im(z):z\in K\}$, $R_{2}=$ $\sup\{\Re(z):z\in K\}$,
    $R_{3}=$ $\inf\{\Re(z):z\in K\}$ and $R=$ $\max\{\mid R_{2}\mid,\mid
    R_{3}\mid\}$. Since $K$ is compact $\varepsilon\neq 0$, $R_{1}<+\infty$ and also $R<+\infty$.
    Let $C$ be the center of the circle passing through
    the points $\frac{\varepsilon}{2}i$, $-R-R_{1}+i\varepsilon$ and
    $R+R_{1}+i\varepsilon$. Then $C$ will be on the imaginary axis,
    hence $C=$ $\alpha i$ for some $\alpha\in$ $\mathbb{R}^{+}$ and this
    $\alpha$ satisfies what we want.
\end{proof}

We formulate and prove our approximation scheme as the following
proposition.

\begin{proposition} Let
    $\varphi:\mathbb{H}\rightarrow$ $\mathbb{H}$ be an analytic self-map
    of $\mathbb{H}$ such that $$\varphi(z)=pz+\psi(z),$$ $p > 0$ and
    $\psi\in$ $H^{\infty}$ is such that $\Im(\psi(z))>\epsilon>0$ for
    all $z\in\mathbb{H}$. Then there is an $\alpha\in$ $\mathbb{R}^{+}$
    such that for $C_{\varphi}:H^{2}\rightarrow$ $H^{2}$ we have
\begin{equation*}
    C_{\varphi} = V_{p}\sum_{n=0}^{\infty}T_{\tau^{n}}D_{\vartheta_{n}},
\end{equation*}
    where the convergence of the series is in operator norm,
    $T_{\tau^{n}}$ is the Toeplitz operator with symbol $\tau^{n}$,
\begin{equation*}
    \tau(x)=i\alpha-\tilde{\psi}(x),\quad
    \tilde{\psi}(x)=\psi(\frac{x}{p}),
\end{equation*}
     $V_{p}$ is the dilation
    operator defined in equation (15) and $D_{\vartheta_{n}}$ is the
    Fourier multiplier with $\vartheta_{n}(t)=$
    $\frac{(-it)^{n}e^{-\alpha t}}{n!}$. \label{prop3}
\end{proposition}

\begin{proof}
    Since for $\varphi(z)=$ $pz+\psi(z)$ where $\psi\in$
    $H^{\infty}$ with $\Im(\psi(z)) > \epsilon > 0$ for all $z\in$
    $\mathbb{H}$ and $p > 0$, we have
\begin{equation*}
    \Im(\varphi^{*}(x))\geq\epsilon>0\quad\textrm{for almost every}\quad x\in\mathbb{R}.
\end{equation*}
    We can use Proposition 3
    for $C_{\varphi}:H^{2}\rightarrow$ $H^{2}$ to have
    \begin{equation*}
    (C_{\varphi}f)(x) = \frac{1}{2\pi i}\int_{-\infty}^{\infty}\frac{f(w)dw}{w-\varphi(x)}=\frac{1}{2\pi i}\int_{-\infty}^{\infty}\frac{f(w)dw}{w-px-\psi(x)}.
    \end{equation*}
    Without loss of generality, we take $p=1$, since if $p\neq 1$ then we have
    \begin{equation}
        (V_{\frac{1}{p}}C_{\varphi})(f)(x) = \frac{1}{2\pi
        i}\int_{-
        \infty}^{\infty}\frac{f(w)dw}{w-x-\tilde{\psi}(x)},
    \end{equation}
     where $\tilde{\psi}(x)=$ $\psi(\frac{x}{p})$ and
    $V_{\beta}f(z)=$ $f(\beta z)$ ($\beta > 0$) is the dilation
    operator. We observe that
    \begin{eqnarray}
        & &\frac{-1}{x-w+\psi(x)} = \frac{-1}{x-w+i\alpha-(i\alpha-\psi(x))}= \\
        &\nonumber &\frac{-1}{(x-w+i\alpha)\left(1-\left(\dfrac{i\alpha-\psi(x)}{x-w+i\alpha}\right)\right)} .
    \end{eqnarray}
    Since $\Im(\psi(z))>\epsilon>0$ for all $z\in\mathbb{H}$ and
    $\psi\in H^{\infty}$, we have $\overline{\psi(\mathbb{H})}$ is
    compact in $\mathbb{H}$, and then by Lemma 4 there is an $\alpha > 0$ such that
\begin{equation*}
    \left\lvert\frac{i\alpha-\psi(x)}{x-w+i\alpha}\right\rvert<\delta<1
\end{equation*}
    for all $x,w\in \mathbb{R}$, so we have
    \begin{equation*}
        \frac{1}{1-\left(\dfrac{i\alpha-\psi(x)}{x-w+i\alpha}\right)} =
        \sum_{n=0}^{\infty}\left(\frac{i\alpha-\psi(x)}{x-w+i\alpha}\right)^{n}.
    \end{equation*}
    Inserting this into equation (21) and then into equation (20), we
    have
    \begin{equation*}
        (C_{\varphi}f)(x) = \sum_{n=0}^{M-1}T_{\tau^{n}}K_{n}f(x) +
        R_{M}f(x) ,
    \end{equation*}
    where $T_{\tau^{n}}f(x)=$ $\tau^{n}(x)f(x)$, $\tau(x)=$ $i\alpha
    -\psi(x)$, $K_{n}$ is as in equation (13) and
\begin{equation*}
    R_{M}f(x)=\frac{1}{2\pi i}T_{\tau^{M+1}}\int_{-\infty}^{\infty}\frac{f(w)dw}{(x-w+i\alpha)^{M}(w-x-\psi(x))}.
\end{equation*}
    By equation (14) we have
\begin{equation*}
    K_{n}f(x)=D_{\vartheta_{n}}f(x)\quad\textrm{and}\quad\vartheta_{n}(t)=\frac{(-it)^{n}e^{-\alpha t}}{n!}.
\end{equation*}
Since $C_{\varphi}$ is bounded it is not difficult to see that
\begin{equation*}
    \parallel R_{M}\parallel\leq\parallel T_{\tau}\parallel\parallel
    C_{\varphi}\parallel\delta^{M}
\end{equation*}
     which implies that $\parallel
    R_{M}\parallel\rightarrow 0$ as $M\rightarrow$ $\infty$. Hence we
    have
\begin{equation*}
    C_{\varphi} = \sum_{n=0}^{\infty}T_{\tau^{n}}D_{\vartheta_{n}} ,
\end{equation*}
    where the convergence is in operator norm.
\end{proof}

\section{a $\Psi$-c*-algebra of operators on hardy spaces of analytic functions}

In the preceding section we have shown that ``quasi-parabolic''
composition operators on the upper half-plane lie in the C*-algebra
generated by certain Toeplitz operators and Fourier multipliers. In
this section we will identify the maximal ideal space of the
C*-algebra generated by Toeplitz operators with a class of symbols
and Fourier multipliers. The C*-algebras generated by multiplication
operators and Fourier multipliers on $L^{2}(\mathbb{R})$ are called
``Pseudo-Differential C*-Algebras'' and they have been studied in a
series of papers by Power (see [22], [23]) and by Cordes and Herman
(see [7]). Our C*-algebra is an analogue of ``Pseudo-differential
C*-algebras'' introduced in [22] and [23]; however our C*-algebra
acts on $H^{2}$ instead of $L^{2}$. Our ``$\Psi$-C*-Algebra'' will
be denoted by $\Psi(A,C([0,\infty]))$ and is defined as
\begin{equation*}
\Psi(A,C([0,\infty]))= C^{*}(\mathcal{T}(A)\cup F),
\end{equation*}
where $A\subseteq L^{\infty}(\mathbb{R})$ is a closed subalgebra of
$L^{\infty}(\mathbb{R})$ and $F$ is as defined by equation (12).

We will now show that if $a\in QC(\mathbb{R})$ and $\vartheta\in
C([0,\infty])$, the commutator $[T_{a},D_{\vartheta}]$ is compact on
$H^{2}(\mathbb{H})$. But before that, we state the following fact
from [21, p. 215] which implies that
\begin{equation*}
PM_{a}-M_{a}P\in K(L^{2})
\end{equation*}
for all $a\in QC$, where $P$ denotes the orthogonal projection of
$L^{2}$ onto $H^{2}$:

\begin{lemma}
    Let $a\in$ $L^{\infty}(\mathbb{T})$ and $P$ be
    the orthogonal projection of $L^{2}(\mathbb{T})$ onto $H^{2}(\mathbb{D})$ then the
    commutator $[P,M_{a}]=$ $PM_{a}-M_{a}P$ is compact
    on $L^{2}(\mathbb{T})$ if and only if $a\in$ $QC$.
    \label{lem4}
\end{lemma}

The following lemma and its proof is a slight modification of Lemma
2.0.15 of [29].

\begin{lemma}
    Let $a\in$ $QC(\mathbb{R})$ and $\vartheta\in$
    $C([0,\infty])$. Then we have
\begin{equation*}
    [T_{a},D_{\vartheta}]=T_{a}D_{\vartheta}-D_{\vartheta}T_{a}\in
    K(H^{2}(\mathbb{H})).
\end{equation*}
     \label{lem3}
\end{lemma}
\begin{proof}
    Let $\tilde{P}:L^{2}(\mathbb{R})\rightarrow H^{2}(\mathbb{H})$ be
    the orthogonal projection of $L^{2}$ onto $H^{2}$ and let $a\in$
    $QC(\mathbb{R})$. Observe that
\begin{equation*}
D_{\chi_{[0,\infty)}}=\tilde{P}
\end{equation*}
    where $\chi_{[0,\infty)}$ is the
    characteristic function of $[0,\infty)$. Let
    $P:L^{2}(\mathbb{T})\rightarrow H^{2}(\mathbb{D})$ be the orthogonal
    projection of $L^{2}$ onto $H^{2}$ on the unit disc. By Lemma 6 and
    by the use of $\Phi$ defined as in equation (4) (observe that $\Phi$
    extends to be an isometric isomorphism from $L^{2}(\mathbb{T})$ onto
    $L^{2}(\mathbb{R})$) we have
\begin{equation*}
    [M_{a},D_{\chi_{[0,\infty)}}]=[M_{a},\tilde{P}]=\Phi\circ [M_{a\circ\mathfrak{C}^{-1}},P]\circ\Phi^{-1}\in K(L^{2}(\mathbb{R})).
\end{equation*}
    Consider $D_{\chi_{[t,\infty)}}$ for $t > 0$ on $L^{2}$:
\begin{eqnarray*}
  & &D_{\chi_{[t,\infty)}}=\mathcal{F}^{-1}M_{\chi_{[t,\infty)}}\mathcal{F}=\mathcal{F}^{-1}S_{-t}M_{\chi_{[0,\infty)}}S_{t}\mathcal{F}=\\
  & &\mathcal{F}^{-1}S_{-t}\mathcal{F}\mathcal{F}^{-1}M_{\chi_{[0,\infty)}}\mathcal{F}\mathcal{F}^{-1}S_{t}\mathcal{F} = M_{e^{-itw}}D_{\chi_{[0,\infty)}}M_{e^{itw}},
\end{eqnarray*}
    where $S_{t}:L^{2}\rightarrow$ $L^{2}$ is the translation operator
\begin{equation*}
    S_{t}f(x)=f(x+t).
\end{equation*}
     Hence we have
\begin{equation*}
        [M_{a},D_{\chi_{[t,\infty)}}]=M_{e^{-itw}}[M_{a},D_{\chi_{[0,\infty)}}]M_{e^{itw}}\in K(L^{2}(\mathbb{R})).
\end{equation*}
    Since the algebra of compact operators is an ideal. So we
    have
\begin{equation*}
    [T_{a},D_{\chi_{[t,\infty)}}]=\tilde{P}[M_{a},D_{\chi_{[t,\infty)}}]|_{H^{2}}\in K(H^{2}(\mathbb{H})).
\end{equation*}
    Consider the characteristic function $\chi_{[t,r)}$ of
    some interval $[t,r)$ where $0 < t < r$. Since
\begin{equation*}
\chi_{[t,r)}=\chi_{[t,\infty)}-\chi_{[r,\infty)}
\end{equation*}
    we have
\begin{equation*}
D_{\chi_{[t,r)}}=D_{\chi_{[t,\infty)}}-D_{\chi_{[r,\infty)}}.
\end{equation*}
    So
\begin{equation*}
[T_{a},D_{\chi_{[t,r)}}]=[T_{a},D_{\chi_{[t,\infty)}}]-[T_{a},D_{\chi_{[r,\infty)}}]\in
K(H^{2}(\mathbb{H})).
\end{equation*}
    Let $\vartheta\in$ $C([0,\infty])$ then for all
    $\varepsilon > 0$ there are $t_{0}= 0 < t_{1}<$...$<t_{n}\in$
    $\mathbb{R}^{+}$ and $c_{1}$,$c_{2}$,...,$c_{n}$, $c_{n+1}\in$
    $\mathbb{C}$ such that
\begin{equation*}
    \|\vartheta-(\sum_{j=1}^{n}c_{j}\chi_{[t_{j-1},t_{j})})-c_{n+1}\chi_{[t_{n},\infty)}\|_{\infty}<\frac{\varepsilon}{2\|
        T_{a}\|}.
\end{equation*}
    Hence we have
        \[
            \|
            [T_{a},D_{\vartheta}]-[T_{a},\sum_{j=1}^{n}c_{j}D_{\chi_{[t_{j-1},t_{j})}}+c_{n+1}D_{\chi_{[t_{n+1},\infty)}}]\|=
        \]
        $$\| [T_{a},D_{\vartheta-(\sum_{j=1}^{n}c_{j}\chi_{[t_{j-1},t_{j}]})-c_{n+1}\chi_{[t_{n},\infty)}}]\|\leq 2\| T_{a}\|\frac{\varepsilon}{2\| T_{a}\|}=\varepsilon$$
    Since
\begin{equation*}
[T_{a},\sum_{j=1}^{n}c_{j}D_{\chi_{[t_{j-1},t_{j})}}+c_{n+1}D_{\chi_{[t_{n+1},\infty)}}]\in
K(H^{2}(\mathbb{H})),
\end{equation*}
    letting $\varepsilon\rightarrow 0$ we have $[T_{a},D_{\vartheta}]\in$ $K(H^{2}(\mathbb{H}))$.
\end{proof}

Now consider the C*-algebra $\Psi(QC(\mathbb{R}),C([0,\infty]))$. By
Douglas' Theorem and Lemma 7, the commutator ideal
$com(\Psi(QC(\mathbb{R}),C([0,\infty]))$ falls inside the ideal of
compact operators $K(H^{2}(\mathbb{H}))$. Since
$\mathcal{T}(C(\dot{\mathbb{R}}))\subset\Psi(QC(\mathbb{R}),C([0,\infty]))$
as in equation (9) we conclude that
\begin{equation*}
    com(\Psi(QC(\mathbb{R}),C([0,\infty])))= K(H^{2}(\mathbb{H})).
\end{equation*}
Therefore we have
\begin{eqnarray}
    \Psi(QC(\mathbb{R}),C([0,\infty]))/K(H^{2}(\mathbb{H}))=\\
    \Psi(QC(\mathbb{R}),C([0,\infty]))/com(\Psi(QC(\mathbb{R}),C([0,\infty]))\nonumber
\end{eqnarray}
and $\Psi(QC(\mathbb{R}),C([0,\infty]))/K(H^{2}(\mathbb{H}))$ is a
commutative C*-algebra with identity. So it is natural to ask for
its maximal ideal space and its Gelfand transform. We will use the
following theorem of Power (see [23]) to characterize its maximal
ideal space:

\begin{theorem} [\scshape Power's Theorem]
    Let $C_{1}$, $C_{2}$ be two C*-subalgebras of $B(H)$ with identity,
    where $H$ is a separable Hilbert space, such that $M(C_{i})\neq$
    $\emptyset$, where $M(C_{i})$ is the space of multiplicative linear
    functionals of $C_{i}$, $i= 1,\,2$ and let $C$ be the C*-algebra that
    they generate. Then for the commutative C*-algebra $\tilde{C}=$
    $C/com(C)$ we have $M(\tilde{C})=$ $P(C_{1},C_{2})\subset$
    $M(C_{1})\times M(C_{2})$, where $P(C_{1},C_{2})$ is defined to be
    the set of points $(x_{1},x_{2})\in$ $M(C_{1})\times M(C_{2})$
    satisfying the condition: \\
    \quad Given $0\leq a_{1} \leq 1$, $0 \leq a_{2} \leq 1$, $a_{1}\in C_{1}$, $a_{2}\in C_{2}$
\begin{equation*}
    x_{i}(a_{i})=1\quad\textrm{with}\quad
        i=1,2\quad\Rightarrow\quad\| a_{1}a_{2}\|=1.
\end{equation*}
    \label{thmpower}
\end{theorem}

Proof of this theorem can be found in [23]. Using Power's
theorem we prove the following result.

\begin{theorem}Let
\begin{equation*}
\Psi(QC(\mathbb{R}),C([0,\infty]))=
C^{*}(\mathcal{T}(QC(\mathbb{R}))\cup F).
\end{equation*}
    Then the C*-algebra
    $\Psi(QC(\mathbb{R}),C([0,\infty]))/K(H^{2}(\mathbb{H}))$ is a
    commutative C*-algebra and its maximal ideal space is
\begin{equation*}
M(\Psi(QC(\mathbb{R}),C([0,\infty])))\cong
(M_{\infty}(QC(\mathbb{R}))\times[0,\infty])\cup(M(QC(\mathbb{R}))\times\{\infty\}),
\end{equation*}
    where
\begin{equation*}
    M_{\infty}(QC(\mathbb{R}))=\{x\in
    M(QC(\mathbb{R})):x|_{C(\dot{\mathbb{R}})}=\delta_{\infty}\quad\textrm{with}\quad\delta_{\infty}(f)=\lim_{t\rightarrow\infty}f(t)\}
\end{equation*}
    is the fiber of $M(QC(\mathbb{R}))$ at $\infty$. \label{thm33}
\end{theorem}

\begin{proof}
    By equation (22) we already know that
    $\Psi(QC(\mathbb{R}),C([0,\infty]))/K(H^{2}(\mathbb{H}))$ is a
    commutative C*-algebra. Since any $x\in M(A)$ vanishes on $com(A)$
    we have
\begin{equation*}
M(A)=M(A/com(A)).
\end{equation*}
    By equation (8)
\begin{equation*}
    \mathcal{T}(QC(\mathbb{R}))/com(\mathcal{T}(QC(\mathbb{R})))=\mathcal{T}(QC(\mathbb{R}))/K(H^{2}(\mathbb{H}))
\end{equation*}
    is isometrically *-isomorphic to $QC(\mathbb{R})$, hence we have
\begin{equation*}
    M(\mathcal{T}(QC(\mathbb{R})))=M(QC(\mathbb{R})).
\end{equation*}
    Now we are ready to use Power's theorem. In our case, $$H=H^{2},
    C_{1}=\mathcal{T}(QC(\mathbb{R})),C_{2}=F\quad\textrm{and}\quad\tilde{C}=\Psi(QC(\mathbb{R}),C([0,\infty]))/K(H^{2}(\mathbb{H})).$$
    We have
\begin{equation*}
    M(C_{1})=M(QC(\mathbb{R}))\quad\textrm{and}\quad M(C_{2})=[0,\infty].
\end{equation*}
    So we need to determine $(x,y)\in$ $M(QC(\mathbb{R}))\times [0,\infty]$ so
    that for all $a\in$ $QC(\mathbb{R})$ and $\vartheta\in$
    $C([0,\infty])$ with $0 < a, \vartheta\leq 1$, we have
\begin{equation*}
    \hat{a}(x)=\vartheta(y)=1\Rightarrow\parallel
    T_{a}D_{\vartheta}\parallel=1\quad\textrm{or}\quad\parallel
    D_{\vartheta}T_{a}\parallel=1.
\end{equation*}
    For any $x\in$ $M(QC(\mathbb{R}))$ consider $\tilde{x}=$
    $x|_{C(\dot{\mathbb{R}})}$ then $\tilde{x}\in$
    $M(C(\dot{\mathbb{R}}))=$ $\dot{\mathbb{R}}$. Hence
    $M(QC(\mathbb{R}))$ is fibered over $\dot{\mathbb{R}}$, i.e.
\begin{equation*}
    M(QC(\mathbb{R}))=\bigcup_{t\in\dot{\mathbb{R}}}M_{t},
\end{equation*}
    where
\begin{equation*}
    M_{t}=\{x\in M(QC(\mathbb{R})):\tilde{x}=x|_{C(\dot{\mathbb{R}})}=\delta_{t}\}.
\end{equation*}
     Let $x\in M(QC(\mathbb{R}))$ such that $x\in M_{t}$ with $t\neq \infty$
    and $y\in$ $[0,\infty)$. Choose $a\in$ $C(\dot{\mathbb{R}})$ and
    $\vartheta\in$ $C([0,\infty])$ such that
\begin{equation*}
    \hat{a}(x)=\phi(t)=\vartheta(y)=1,\quad 0\leq a\leq 1,\quad
    0\leq\vartheta\leq 1,\quad a(z)<1
\end{equation*}
    for all $z\in$
    $\mathbb{R}\backslash\{t\}$ and\quad $\vartheta(w) < 1$\quad for
    all $w\in$ $[0,\infty]\backslash\{y\}$\quad, where both $a$ and
    $\vartheta$ have compact supports. Consider $\parallel
    T_{a}D_{\vartheta}\parallel_{H^{2}}$. Let $\tilde{\vartheta}$ be
\begin{equation*}
    \tilde{\vartheta}(w)=
        \begin{cases}
            \vartheta(w)\qquad\textrm{if}\qquad w\geq 0\\
            \vartheta(-w)\qquad\textrm{if}\qquad w < 0
        \end{cases}
\end{equation*}
    then
\begin{equation*}
    PM_{a}D_{\tilde{\vartheta}}|_{H^{2}}=T_{a}D_{\vartheta},
\end{equation*}
    where $P:L^{2}\rightarrow$ $H^{2}$ is the orthogonal projection of
    $L^{2}$ onto $H^{2}$. So we have
\begin{equation*}
    \parallel T_{a}D_{\vartheta}\parallel_{H^{2}}\leq\parallel M_{a}D_{\tilde{\vartheta}}\parallel_{L^{2}}.
\end{equation*}
     By a result of
    Power (see [22] and also [29]) under these conditions we have
    \begin{equation}
        \parallel M_{a}D_{\tilde{\vartheta}}\parallel_{L^{2}}<
        1\Rightarrow\parallel T_{a}D_{\vartheta}\parallel_{H^{2}}<
        1\Rightarrow (x,y)\not\in M(\tilde{C}),
    \end{equation}

    so if $(x,y)\in M(\tilde{C})$, then either $y= \infty$ or $x\in
    M_{\infty}(QC(\mathbb{R}))$.

    Let $y= \infty$ and $x\in M(QC(\mathbb{R}))$. Let $a\in
    QC(\mathbb{R})$ and $\vartheta\in C([0,\infty])$ such that
\begin{equation*}
    0\leq a,\vartheta\leq 1\quad\textrm{and}\quad\hat{a}(x)=
    \vartheta(y)= 1.
\end{equation*}
    Consider
    \begin{eqnarray}
        \parallel
        D_{\vartheta}T_{a}\parallel_{H^{2}}=\parallel\mathcal{F}D_{\vartheta}T_{a}\mathcal{F}^{-1}\parallel_{L^{2}([0,\infty))}=\parallel
        M_{\vartheta}\mathcal{F}T_{a}\mathcal{F}^{-1}\parallel_{L^{2}([0,\infty))}\\
        =\parallel
        M_{\vartheta}\mathcal{F}(\mathcal{F}^{-1}M_{\chi_{[0,\infty)}}\mathcal{F})M_{a}\mathcal{F}^{-1}\parallel_{L^{2}([0,\infty))}=\parallel
        M_{\vartheta}\mathcal{F}M_{a}\mathcal{F}^{-1}\parallel_{L^{2}([0,\infty))} . \nonumber
    \end{eqnarray}

    Choose $f\in$ $L^{2}([0,\infty))$ with $\parallel
    f\parallel_{L^{2}([0,\infty))}=1$ such that
\begin{equation*}
    \parallel(\mathcal{F}M_{a}\mathcal{F}^{-1})f\parallel\geq
    1-\varepsilon
\end{equation*}
    for given $\varepsilon > 0$. Since
    $\vartheta(\infty)= 1$ there exists $w_{0} > 0$ so that
\begin{equation*}
    1-\varepsilon\leq\vartheta(w)\leq 1\quad\forall w\geq w_{0}.
\end{equation*}
    Let $t_{0}\geq w_{0}$. Since the support of
    $(S_{-t_{0}}\mathcal{F}M_{a}\mathcal{F}^{-1})f$ lies in
    $[t_{0},\infty)$ where $S_{t}$ is the translation by $t$, we have
    \begin{eqnarray}
    & &\parallel M_{\vartheta}(S_{-t_{0}}\mathcal{F}M_{a}\mathcal{F}^{-1})f\parallel_{2}\\
    \geq&\nonumber &\inf\{\vartheta(w): w\in (w_{0},\infty)\}\parallel(\mathcal{F}M_{a}\mathcal{F}^{-1})f\parallel_{2}\geq(1-\varepsilon)^{2}
    \end{eqnarray}
    Since
\begin{equation*}
    S_{-t_{0}}\mathcal{F}M_{a}\mathcal{F}^{-1}=\mathcal{F}M_{a}\mathcal{F}^{-1}S_{-t_{0}}
\end{equation*}
    and $S_{-t_{0}}$ is an isometry on $L^{2}([0,\infty))$ by equations
    (24) and (25), we conclude that
\begin{equation*}
    \parallel M_{\vartheta}\mathcal{F}M_{a}\mathcal{F}^{-1}\parallel_{L^{2}([0,\infty))}=\parallel
    D_{\vartheta}T_{a}\parallel_{H^{2}}= 1\Rightarrow (x,\infty)\in
    M(\tilde{C})\quad\forall x\in M(C_{1}).
\end{equation*}
    Now let $x\in M_{\infty}(QC(\mathbb{R}))$ and $y\in [0,\infty]$. Let
    $a\in QC(\mathbb{R})$ and $\vartheta\in C([0,\infty])$ such that
\begin{equation*}
    \hat{a}(x)=\vartheta(y)=1\quad\textrm{and}\quad
    0\leq\phi,\vartheta\leq 1.
\end{equation*}
    By a result of Sarason (see [27]
    lemmas 5 and 7) for a given $\varepsilon > 0$ there is a $\delta> 0$ so that
    \begin{equation}
    \mid\hat{a}(x)-\frac{1}{2\delta}\int_{-\delta}^{\delta}a\circ\mathfrak{C}^{-1}(e^{i\theta})d\theta\mid\leq\varepsilon.
    \end{equation}
    Since $\hat{a}(x)= 1$ and $0\leq a\leq 1$, this implies
    that for all $\varepsilon > 0$ there exists $w_{0} > 0$ such that
    $1-\varepsilon\leq$ $a(w)\leq 1$ for a.e.\quad$w$ with $\mid
    w\mid > w_{0}$. Let $\tilde{\vartheta}$ be
\begin{equation*}
\tilde{\vartheta}(w)=
    \begin{cases}
    \vartheta(w)\qquad\textrm{if}\qquad w\geq 0\\
    0\qquad\textrm{if}\qquad w < 0
    \end{cases} .
\end{equation*}
 Then we have
\begin{equation*}
D_{\vartheta}T_{a}=D_{\tilde{\vartheta}}M_{a}.
\end{equation*}
    Let $\varepsilon > 0$ be given. Let
    $g\in$ $H^{2}$ so that $\parallel g\parallel_{2}= 1$ and
    $\parallel D_{\tilde{\vartheta}}g\parallel_{2}\geq 1-\varepsilon$.
    Since $1-\varepsilon\leq a(w)\leq 1$ for a.e. $w$ with
    $\mid w\mid > w_{0}$, there is a $w_{1} > 2w_{0}$ so that
\begin{equation*}
\parallel S_{w_{1}}g-M_{a}S_{w_{1}}g\parallel_{2}\leq 2\varepsilon.
\end{equation*}
    We have $\parallel D_{\tilde{\vartheta}}\parallel = 1$ and this implies that
    \begin{equation}
        \parallel D_{\tilde{\vartheta}}S_{w_{1}}g-D_{\tilde{\vartheta}}M_{a}S_{w_{1}}g\parallel_{2}\leq2\varepsilon .
    \end{equation}
    Since $S_{w}D_{\tilde{\vartheta}}= D_{\tilde{\vartheta}}S_{w}$ and
    $S_{w}$ is unitary for all $w\in \mathbb{R}$, we have
    $$\parallel D_{\tilde{\vartheta}}M_{a}S_{w_{1}}g\parallel_{2}\geq 1-3\varepsilon$$
    and $(x,y)\in M(\tilde{C})$ for all $x\in M_{\infty}(C_{1})$.
\end{proof}

 The Gelfand transform $\Gamma$ of
$\Psi(QC(\mathbb{R}),C([0,\infty]))/K(H^{2}(\mathbb{H}))$ looks like
\begin{equation}
\Gamma\left(\left[\sum_{j=1}^{\infty}T_{a_{j}}D_{\vartheta_{j}}\right]\right)(x,t)=
\begin{cases}
\sum_{j=1}^{\infty}\hat{a_{j}}(x)\hat{\vartheta_{j}}(t)\quad\textrm{if}\quad
x\in
M_{\infty}(QC(\mathbb{R})) \\
\sum_{j=1}^{\infty}\hat{a_{j}}(x)\hat{\vartheta_{j}}(\infty)\quad\textrm{if}\quad
t=\infty
\end{cases}
\end{equation}

\section{main results}

In this section we characterize the essential spectra of
quasi-parabolic composition operators with translation functions in
$QC$ class which is the main aim of the paper. In doing this we will
heavily use Banach algebraic methods. We start with the following
proposition from Hoffman's book (see [15] pp.171):

\begin{proposition} Let $f$ be a function in
$A\subseteq L^{\infty}(\mathbb{T})$ where $A$ is a closed
*-subalgebra of $L^{\infty}(\mathbb{T})$ which contains
$C(\mathbb{T})$. The range of $\hat{f}$ on the fiber $M_{\alpha}(A)$
consists of all complex numbers $\zeta$ with this property: for each
neighborhood $N$ of $\alpha$ and each $\varepsilon > 0$, the set
$$\{\mid f-\zeta\mid<\varepsilon\}\cap N$$ has positive Lebesgue
measure. \label{prop4}
\end{proposition}
Hoffman states and proves Proposition 10 for
$A=L^{\infty}(\mathbb{T})$ but in fact his proof works for a general
C*-subalgebra of $L^{\infty}(\mathbb{T})$ that contains
$C(\mathbb{T})$.

Using a result of Shapiro [30] we deduce the following lemma that
might be regarded as the upper half-plane version of that result:
\begin{lemma}
If $\psi\in$ $QC(\mathbb{R})\cap H^{\infty}(\mathbb{H})$ we have
$$\mathcal{R}_{\infty}(\psi)=\mathcal{C}_{\infty}(\psi)$$ where $\mathcal{C}_{\infty}(\psi)$ is the cluster set of $\psi$ at
infinity which is defined as the set of points $z\in$ $\mathbb{C}$
for which there is a sequence $\{z_{n}\}\subset$ $\mathbb{H}$ so
that $z_{n}\rightarrow$ $\infty$ and $\psi(z_{n})\rightarrow$ $z$.
\end{lemma}
\begin{proof}
 Since the pullback measure $\lambda_{0}(E)=$ $\mid\mathfrak{C}(E)\mid$ is
absolutely continuous with respect to the Lebesgue measure $\lambda$
on $\mathbb{R}$ where $\mid\cdot\mid$ denotes the Lebesgue measure
on $\mathbb{T}$ and $E$ denotes a Borel subset of $\mathbb{R}$, we
have if $\psi\in$ $L^{\infty}(\mathbb{R})$ then
\begin{equation}
\mathcal{R}_{\infty}(\psi)=\mathcal{R}_{1}(\psi\circ
\mathfrak{C}^{-1}).
\end{equation}
By the result of Shapiro (see[30]) if $\psi\in QC(\mathbb{R})\cap
H^{\infty}(\mathbb{H})$ then we have
\begin{equation*}
\mathcal{R}_{1}(\psi\circ
\mathfrak{C}^{-1})=\mathcal{C}_{1}(\psi\circ\mathfrak{C}^{-1})
\end{equation*}
Since
\begin{equation*}
\mathcal{C}_{1}(\psi\circ\mathfrak{C}^{-1})=\mathcal{C}_{\infty}(\psi)
\end{equation*}
we have
\begin{equation*}
\mathcal{R}_{\infty}(\psi)=\mathcal{C}_{\infty}(\psi)
\end{equation*}
\end{proof}

Firstly we have the following result on the upper half-plane:

\begin{thma}
Let $\psi\in$ $QC(\mathbb{R})\cap H^{\infty}(\mathbb{H})$ such that
$\Im(\psi(z)) > \epsilon > 0$ for all $z\in$ $\mathbb{H}$ then for
$\varphi(z)=$ $z+\psi(z)$ we have
\begin{itemize}
\item (i)\quad $C_{\varphi}:$ $H^{2}(\mathbb{H})\rightarrow$
$H^{2}(\mathbb{H})$ is essentially normal \\
\item (ii)\quad$\sigma_{e}(C_{\varphi})=\{e^{izt}:t\in [0,\infty],
z\in\mathcal{C}_{\infty}(\psi)=\mathcal{R}_{\infty}(\psi)\}\cup\{0\}$
\end{itemize}
where\quad $\mathcal{C}_{\infty}(\psi)$ and
$\mathcal{R}_{\infty}(\psi)$ are the set of cluster points and the
local essential range of $\psi$ at $\infty$ respectively.
\end{thma}

\begin{proof} By Proposition 5 we have the following series
expansion for $C_{\varphi}$:
\begin{equation}
C_{\varphi}=\sum_{j=0}^{\infty}\frac{1}{j!}(T_{\tau})^{j}D_{(-it)^{j}e^{-\alpha
t}}
\end{equation}
where $\tau(z)=$ $i\alpha-\psi(z)$. So we conclude that if $\psi\in$
$QC(\mathbb{R})\cap H^{\infty}(\mathbb{H})$ with $\Im({\psi(z)})>$
$\epsilon > 0$ then
\begin{equation*}
C_{\varphi}\in\Psi(QC(\mathbb{R}),C([0,\infty]))
\end{equation*}
where $\varphi(z)=$ $z+\psi(z)$. Since
$\Psi(QC(\mathbb{R}),C([0,\infty]))/K(H^{2}(\mathbb{H}))$ is
commutative, for any $T\in\Psi(QC(\mathbb{R}),C([0,\infty]))$ we
have $T^{*}\in \Psi(QC(\mathbb{R}),C([0,\infty]))$ and
\begin{equation}
[TT^{*}]=[T][T^{*}]=[T^{*}][T]=[T^{*}T].
\end{equation}
This implies that $(TT^{*}-T^{*}T)\in K(H^{2}(\mathbb{H}))$. Since
$C_{\varphi}\in\Psi(QC(\mathbb{R}),C([0,\infty]))$ we also have
\begin{equation*}
(C_{\varphi}^{*}C_{\varphi}-C_{\varphi}C_{\varphi}^{*})\in
K(H^{2}(\mathbb{H})).
\end{equation*}
This proves (i).

For (ii) we look at the values of $\Gamma[C_{\varphi}]$ at
$M(\Psi(QC(\mathbb{R}),C([0,\infty]))/K(H^{2}(\mathbb{H})))$ where
$\Gamma$ is the Gelfand transform of
$\Psi(QC(\mathbb{R}),C([0,\infty]))/K(H^{2}(\mathbb{H}))$. By
Theorem 9 we have
\begin{equation*}
M(\Psi(QC(\mathbb{R}),C([0,\infty]))/K(H^{2}(\mathbb{H})))=(M(QC(\mathbb{R}))\times\{\infty\})\cup(M_{\infty}(QC(\mathbb{R}))\times
[0,\infty]).
\end{equation*}
By equations (28) and (30) we have the Gelfand transform
$\Gamma[C_{\varphi}]$ of $C_{\varphi}$ at $t=$ $\infty$ as
\begin{equation}
(\Gamma[C_{\varphi}])(x,\infty)=\sum_{j=0}^{\infty}\frac{1}{j!}\hat{\tau}(x)\vartheta_{j}(\infty)=0\quad\forall
x\in M(QC(\mathbb{R}))
\end{equation}
since $\vartheta_{j}(\infty)=0$ for all $j\in\mathbb{N}$ where
$\vartheta_{j}(t)=(-it)^{j}e^{-\alpha t}$. We calculate
$\Gamma[C_{\varphi}]$ of $C_{\varphi}$ for $x\in$
$M_{\infty}(QC(\mathbb{R}))$ as
%\begin{align*}
%\Gamma(x,t)([C_{\varphi}])&=\Gamma(x,t)([\sum_{j=0}^{\infty}\frac{1}{j!}(T_{\tau})^{j}D_{(-it)^{j}e^{-\alpha
%t}}])\\
%&=\sum_{j=0}^{\infty}\frac{1}{j!}\hat{\tau}(x)^{j}(-it)^{j}e^{-\alpha
%t}&=e^{i\hat{\psi}(x)t}
%\end{align*}
\begin{eqnarray}
& &(\Gamma[C_{\varphi}])(x,t)=\\
 &\nonumber &\left(\Gamma\left[\sum_{j=0}^{\infty}\frac{1}{j!}(T_{\tau})^{j}D_{(-it)^{j}e^{-\alpha
t}}\right]\right)(x,t)=\sum_{j=0}^{\infty}\frac{1}{j!}\hat{\tau}(x)^{j}(-it)^{j}e^{-\alpha
t}=e^{i\hat{\psi}(x)t}
\end{eqnarray}
 for all
$x\in$ $M_{\infty}(QC(\mathbb{R}))$ and $t\in$ $[0,\infty]$. So we
have $\Gamma[C_{\varphi}]$ as the following:

\begin{eqnarray}
& &\Gamma([C_{\varphi}])(x,t)=
\begin{cases}
e^{i\hat{\psi}(x)t}\quad\textrm{if}\quad x\in
M_{\infty}(QC(\mathbb{R})) \\
0\quad\qquad\textrm{if}\quad t=\infty
\end{cases}
\end{eqnarray}
Since $\Psi=\Psi(QC(\mathbb{R}),C([0,\infty]))/K(H^{2}(\mathbb{H}))$
is a commutative Banach algebra with identity, by equations (2) and
(34) we have
\begin{eqnarray}
\sigma_{\Psi}([C_{\varphi}])=\{\Gamma[C_{\varphi}](x,t):(x,t)\in
M(\Psi(QC(\mathbb{R}),C([0,\infty]))/K(H^{2}))\}=\\
\{e^{i x(\psi)t}:x\in M_{\infty}(QC(\mathbb{R})),t\in\nonumber
[0,\infty)\}\cup\{0\}
\end{eqnarray}
Since $\Psi$ is a closed *-subalgebra of the Calkin algebra
$\mathcal{B}(H^{2}(\mathbb{H}))/K(H^{2}(\mathbb{H}))$ which is also
a C*-algebra, by equation (1) we have
\begin{equation}
\sigma_{\Psi}([C_{\varphi}])=\sigma_{\mathcal{B}(H^{2})/K(H^{2})}([C_{\varphi}]).
\end{equation}
But by definition
$\sigma_{\mathcal{B}(H^{2})/K(H^{2})}([C_{\varphi}])$ is the
essential spectrum of $C_{\varphi}$. Hence we have
\begin{equation}
\sigma_{e}(C_{\varphi})=\{e^{i\hat{\psi}(x)t}:x\in
M_{\infty}(QC(\mathbb{R})),t\in [0,\infty)\}\cup\{0\}.
\end{equation}

Now it only remains for us to understand what the set
$\{\hat{\psi}(x)=x(\psi):x\in M_{\infty}(QC(\mathbb{R}))\}$ looks
like, where $M_{\infty}(QC(\mathbb{R}))$ is as defined in Theorem 9.
By Proposition 10 and equation (29) we have
\begin{equation*}
\{\hat{\psi}(x):x\in
M_{\infty}(QC(\mathbb{R}))\}=\{\hat{\psi\circ\mathfrak{C}^{-1}}(x):x\in
M_{1}(QC)\}=\mathcal{R}_{1}(\psi\circ\mathfrak{C}^{-1})=\mathcal{R}_{\infty}(\psi).
\end{equation*}
By Lemma 11 we have
\begin{eqnarray*}
& &\sigma_{e}(C_{\varphi})=\{(\Gamma[C_{\varphi}])(x,t):(x,t)\in
M(\Psi(QC(\mathbb{R}),C([0,\infty]))/K(H^{2}(\mathbb{H})))\}=\\
& &\{e^{izt}:t\in [0,\infty),
z\in\mathcal{C}_{\infty}(\psi)=\mathcal{R}_{\infty}(\psi)\}\cup\{0\}
\end{eqnarray*}
\end{proof}

\begin{thmb}
 Let $\varphi:\mathbb{D}\rightarrow$ $\mathbb{D}$ be an
analytic self-map of $\mathbb{D}$ such that
\begin{equation*}
\varphi(z)=\frac{2iz+\eta(z)(1-z)}{2i+\eta(z)(1-z)}
\end{equation*}
where $\eta\in$ $H^{\infty}(\mathbb{D})$ with $\Im(\eta(z)) >
\epsilon> 0$ for all $z\in$ $\mathbb{D}$. If $\eta\in$ $QC\cap
H^{\infty}$ then we have
\begin{itemize}
\item (i)\quad $C_{\varphi}:$ $H^{2}(\mathbb{D})\rightarrow$
$H^{2}(\mathbb{D})$ is essentially normal
\item (ii)\quad$\sigma_{e}(C_{\varphi})=\{e^{izt}:t\in [0,\infty],
z\in\mathcal{C}_{1}(\eta)=\mathcal{R}_{1}(\eta)\}\cup\{0\}$
\end{itemize}
where\quad $\mathcal{C}_{1}(\eta)$ and $\mathcal{R}_{1}(\eta)$ are
the set of cluster points and the local essential range of $\eta$ at
 $1$ respectively.
\end{thmb}

 \begin{proof}
Using the isometric isomorphism
$\Phi:H^{2}(\mathbb{D})\longrightarrow$ $H^{2}(\mathbb{H})$
introduced in section 2, if $\varphi:\mathbb{D}\rightarrow$
$\mathbb{D}$ is of the form
\begin{equation*}
\varphi(z)=\frac{2iz+\eta(z)(1-z)}{2i+\eta(z)(1-z)}
\end{equation*}
where $\eta\in$ $H^{\infty}(\mathbb{D})$ satisfies $\Im(\eta(z))>$
$\delta
> 0$ then, by equation (9), for $\tilde{\varphi}=$
$\mathfrak{C}^{-1}\circ\varphi\circ\mathfrak{C}$ we have
$\tilde{\varphi}(z)=$ $z+\eta\circ\mathfrak{C}(z)$ and
\begin{equation}
\Phi\circ
C_{\varphi}\circ\Phi^{-1}=C_{\tilde{\varphi}}+T_{\frac{\eta\circ\mathfrak{C}(z)}{z+i}}C_{\tilde{\varphi}}.
\end{equation}
For $\eta\in$ $QC$ we have both
\begin{equation*}
C_{\tilde{\varphi}}\in\Psi(QC(\mathbb{R}),C([0,\infty]))\quad\textrm{and}\quad
T_{\frac{\eta\circ\mathfrak{C}(z)}{z+i}}\in\Psi(QC(\mathbb{R}),C([0,\infty]))
\end{equation*}
and hence
\begin{equation*}
 \Phi\circ
C_{\varphi}\circ\Phi^{-1}\in\Psi(QC(\mathbb{R}),C([0,\infty])).
\end{equation*}
Since $\Psi(QC(\mathbb{R}),C([0,\infty]))/K(H^{2})$ is commutative
and $\Phi$ is an isometric isomorphism, (i) follows from the
argument following equation (5)($C_{\varphi}$ is essentially normal
if and only if $\Phi\circ C_{\varphi}\circ\Phi^{-1}$ is essentially
normal) and by equation (31).

For (ii) we look at the values of $\Gamma[\Phi\circ
C_{\varphi}\circ\Phi^{-1}]$ at
$M(\Psi(QC(\mathbb{R}),C([0,\infty]))/K(H^{2}))$ where $\Gamma$ is
the Gelfand transform of
$\Psi(QC(\mathbb{R}),C([0,\infty]))/K(H^{2})$. Again applying the
Gelfand transform for
\begin{equation*}
(x,\infty)\in M(\Psi(QC(\mathbb{R}),C([0,\infty]))/K(H^{2}))\subset
M(QC(\mathbb{R}))\times [0,\infty]
\end{equation*}
 we have
\begin{eqnarray*}
& &(\Gamma[\Phi\circ C_{\varphi}\circ\Phi^{-1}])(x,\infty)=\\
 & &(\Gamma[C_{\tilde{\varphi}}])(x,\infty)+((\Gamma[T_{\eta\circ\mathfrak{C}}])(x,\infty))((\Gamma[T_{\frac{1}{i+z}}])(x,\infty))((\Gamma[C_{\tilde{\varphi}}])(x,\infty))
\end{eqnarray*}

Appealing to equation (32) we have
$(\Gamma[C_{\tilde{\varphi}}])(x,\infty)= 0$ for all $x\in$
$M(QC(\mathbb{R}))$ hence we have
\begin{equation*}
(\Gamma[\Phi\circ C_{\varphi}\circ\Phi^{-1}])(x,\infty)=0
\end{equation*}
 for all $x\in$ $M(QC(\mathbb{R}))$. Applying the Gelfand transform for
$$(x,t)\in M_{\infty}(QC)\times [0,\infty]\subset
M(\Psi(QC(\mathbb{R}),C([0,\infty]))/K(H^{2})$$ we have
\begin{equation*}
(\Gamma[\Phi\circ
C_{\varphi}\circ\Phi^{-1}])(x,t)=(\Gamma[C_{\tilde{\varphi}}])(x,t)+((\Gamma[T_{\eta\circ\mathfrak{C}}])(x,t))((\Gamma[T_{\frac{1}{i+z}}])(x,t))((\Gamma[C_{\tilde{\varphi}}])(x,t)).
\end{equation*}
Since $x\in$ $M_{\infty}(QC(\mathbb{R}))$ we have
\begin{equation*}
(\Gamma[T_{\frac{1}{i+z}}])(x,t)=\hat{\left(\frac{1}{i+z}\right)}(x)=x\left(\frac{1}{i+z}\right)=0.
\end{equation*}
Hence we have
\begin{equation*}
(\Gamma[\Phi\circ
C_{\varphi}\circ\Phi^{-1}])(x,t)=(\Gamma[C_{\tilde{\varphi}}])(x,t)
\end{equation*}
for all $(x,t)\in$ $M_{\infty}(QC(\mathbb{R}))\times [0,\infty]$.
Moreover we have
\begin{equation*}
(\Gamma[\Phi\circ
C_{\varphi}\circ\Phi^{-1}])(x,t)=(\Gamma[C_{\tilde{\varphi}}])(x,t)
\end{equation*}
for all $(x,t)\in M(\Psi(QC(\mathbb{R}),C([0,\infty]))/K(H^{2}))$.
Therefore by similar arguments in Theorem A (equations (35) and
(36)) we have
\begin{equation*}
\sigma_{e}(\Phi\circ
C_{\varphi}\circ\Phi^{-1})=\sigma_{e}(C_{\tilde{\varphi}}).
\end{equation*}
 By Theorem A (together with equation (37)) and equation (5) we have
\begin{equation*}
\sigma_{e}(C_{\varphi})=\sigma_{e}(C_{\tilde{\varphi}})=\{e^{izt}:z\in\mathcal{R}_{\infty}(\eta\circ\mathfrak{C})=\mathcal{R}_{1}(\eta),t\in
[0,\infty]\}.
\end{equation*}
\end{proof}

We add a few remarks on these theorems:

1. In case $\psi\in$ $H^{\infty}(\mathbb{H})\cap
C(\dot{\mathbb{R}})\subset QC(\mathbb{R})$ we recall that
\begin{equation*}
 \sigma_{e}(C_{\varphi})=\{e^{iz_{0}t}:t\in
[0,\infty),z_{0}=\lim_{x\rightarrow\infty}\psi(x)\}\cup\{0\}
\end{equation*}
 where
$\varphi(z)=$ $z+\psi(z)$ with $\Im(\psi(z)) > \delta > 0$ for all
$z\in$ $\mathbb{H}$. Hence we recapture a result from the work of
Kriete and Moorhouse [17] and also from the work of Bourdon, Levi,
Narayan and Shapiro [3]. However theorems A and B allow us to
compute the essential spectra of operators not considered by [3] or
[17]. To illustrate this point we give an example of $\psi\in$
$QC(\mathbb{R})\cap H^{\infty}(\mathbb{H})$ so that $\psi\not\in$
$C(\dot{\mathbb{R}})$. Recall that an analytic function $f$ on the
unit disc is in the Dirichlet space if and only if
$\int_{\mathbb{D}}\mid f'(z)\mid^{2}dA(z) < \infty$. We will use the
following proposition (see also [31]):

\begin{proposition}\footnote{ The author is indebted to
Professor Joe Cima for pointing out this
proposition.}
Every bounded analytic function of the unit
disc that is in the Dirichlet space is in $QC$
\end{proposition}

\begin{proof}
By the VMO version of Fefferman's theorem on BMO (see Chapter 5 of
[27]): $f\in$ $VMOA$ if and only if $\mid f'(z)\mid^{2}(1-\mid
z\mid)dA(z)$ is a vanishing Carleson measure. And by a result of
Sarason (see [25] and [27]) we have $VMOA\cap L^{\infty}=$
$([H^{\infty}+C(\mathbb{T})]\cap[\overline{H^{\infty}}+C(\mathbb{T})])\cap
H^{\infty}$. So we only need to show that if $f$ is in the Dirichlet
space then $\mid f'(z)\mid^{2}(1-\mid z\mid)dA(z)$ is a vanishing
Carleson measure: Let $S(I)=$ $\{re^{i\theta}:e^{i\theta}\in
I\quad\textrm{and}\quad 1-\frac{\mid I\mid}{2}<r<1\}$ be the
Carleson window associated to the arc $I\subseteq$ $\mathbb{T}$
where $\mid I\mid$ denotes the length of the arc $I$. Then for any
$z\in$ $S(I)$ we have $1-\mid z\mid < \mid I\mid$. So we have
$$\int_{S(I)}\mid f'(z)\mid^{2}(1-\mid
z\mid)dA(z)\leq\mid I\mid\int_{S(I)}\mid f'(z)\mid^{2}dA(z).$$ Let
$I_{n}$ be a sequence of decreasing arcs so that $\mid
I_{n}\mid\rightarrow 0$ and let $\chi_{S(I_{n})}$ be the
characteristic function of $S(I_{n})$. Then $\chi_{S(I_{n})}$ is in
the unit ball of $L^{\infty}(\mathbb{D},dA)$ which is the dual of
$L^{1}(\mathbb{D},dA)$. Hence by Banach-Alaoglu theorem, there is a
subsequence $\chi_{S(I_{n_{k}})}$ so that
\begin{equation*}
\int_{\mathbb{D}}\chi_{S(I_{n_{k}})}g dA=\int_{S(I_{n_{k}})}g
dA\rightarrow\int_{\mathbb{D}}\phi g dA
\end{equation*}
as $k\rightarrow$ $\infty$ for all $g\in$ $L^{1}(\mathbb{D},dA)$ for
some $\phi\in$ $L^{\infty}(\mathbb{D},dA)$. In particular if we take
$g\in$ $L^{\infty}(\mathbb{D},dA)\subset$ $L^{1}(\mathbb{D},dA)$ we
observe that
\begin{equation*}
\mid\int_{S(I_{n_{k}})}g dA\mid\leq A(S(I_{n_{k}}))\| g\|_{\infty}
\end{equation*}
so we have $\mid\int_{S(I_{n_{k}})}g dA\mid\rightarrow 0$ as
$k\rightarrow$ $\infty$ if $g\in$ $L^{\infty}(\mathbb{D},dA)$. Since
$L^{\infty}(\mathbb{D},dA)$ is dense in $L^{1}(\mathbb{D},dA)$ we
have $\int_{\mathbb{D}}\phi g dA= 0$ for all $g\in$
$L^{1}(\mathbb{D},dA)$ and hence $\phi\equiv 0$. Since $\mid
f'\mid^{2}\in$ $L^{1}(\mathbb{D},dA)$ we have $\int_{S(I_{n})}\mid
f'\mid^{2}dA\rightarrow 0$ as $n\rightarrow$ $\infty$. Since
$L^{1}(\mathbb{D},dA)$ is separable the weak* topology of the unit
ball of $L^{\infty}(\mathbb{D},dA)$ is metrizable. This implies that
$$\lim_{\mid I\mid\rightarrow 0}\int_{S(I)}\mid f'\mid^{2}dA=0.$$
This proves the proposition.
\end{proof}

We are ready to construct our example of a ``quasi-parabolic"
composition operator which has thick essential spectrum: Let $D$ be
the simply connected region bounded by the curve
$\alpha:[-\pi,\pi]\rightarrow$ $\mathbb{C}$ such that $\alpha$ is
continuous on $[-\pi,0)\cup$ $(0,\pi]$ with
$$\alpha(t)=
    \begin{cases}
        3i+i(t+i\frac{(b-a)}{2}\sin(\frac{3\pi^{2}}{4t}))+\frac{(a+b)}{2}\quad\textrm{if}\quad
        t\in
        (0,\frac{\pi}{2}]\\
        i(3+t)+\frac{(a+b)}{2}\quad\textrm{if}\quad t\in [-\frac{\pi}{2},0]\\
    \end{cases}$$
and $\alpha(\pi) = \alpha(-\pi) =
(1+\frac{\pi}{2})+i(3+\frac{\pi}{2})+\frac{(a+b)}{2}$, $a < b$. By
the Riemann mapping theorem there is a conformal mapping
$\tilde{\psi}:\mathbb{D} \rightarrow D$ that is bi-holomorphic. One
can choose $\tilde{\psi}$ to satisfy $\lim_{\theta\rightarrow
0^{-}}\tilde{\psi}(e^{i\theta}) = 3i$. Since $D$ has finite area and
$\tilde{\psi}$ is one-to-one and onto, $\tilde{\psi}$ is in the
Dirichlet space and hence, by Proposition 12, $\tilde{\psi} \in QC$.
Let $\psi = \tilde{\psi}\circ\mathfrak{C}$. Then $\psi\in
QC(\mathbb{R})$ and $\psi\not \in H^{\infty}(\mathbb{H})\cap
C(\dot{\mathbb{R}})$. We observe that
\begin{equation*}
\mathcal{C}_{\infty}(\psi)=\{3i+x:x\in [a,b]\}.
\end{equation*}
So for $C_{\varphi}:H^{2}(\mathbb{H}) \rightarrow H^{2}(\mathbb{H})$
and for $C_{\tilde{\varphi}}:H^{2}(\mathbb{D}) \rightarrow
H^{2}(\mathbb{D})$, we have
\begin{equation*}
\sigma_{e}(C_{\varphi}) =
\sigma_{e}(C_{\tilde{\varphi}})=\{e^{i(3i+x)t}:t\in [0,\infty),x\in
[a,b]\}\cup\{0\}
\end{equation*}
where
\begin{equation*}
\varphi(z) = z+\psi(z)\quad\textrm{and}\quad\tilde{\varphi} =
\mathfrak{C}\circ\varphi\circ\mathfrak{C}^{-1}.
\end{equation*}

2. We recall that for any $a$, $\vartheta\in C(\dot{\mathbb{R}})$
the operator $M_{a}D_{\vartheta}$ is Hilbert Schmidt on
$L^{2}(\mathbb{R})$(see pp.482-484 of [22]). Hence for $a\in
C(\dot{\mathbb{R}})$ and $\vartheta\in C([0,\infty])$ with
$a(\infty)=\vartheta(\infty)=0$ the operator $T_{a}D_{\vartheta}$ is
compact on $H^{2}(\mathbb{H})$. Combining this fact with equations
(30) and (38) we conclude that for
$\varphi:\mathbb{D}\rightarrow\mathbb{D}$ with
\begin{equation*}
\tilde{\varphi}(z)=\mathfrak{C}^{-1}\circ\varphi\circ\mathfrak{C}(z)=z+\psi(z),
\end{equation*}
$\psi\in H^{\infty}$ and $\Im(\psi(z))>\epsilon>0$, the operator
$\Phi\circ C_{\varphi}\circ\Phi^{-1}-C_{\tilde{\varphi}}$ is compact
on $H^{2}(\mathbb{H})$.

\section{further results}

In this last part of the paper we will prove a more general result
about $C_{\varphi}$ with $\varphi(z)=$ $z+\psi(z)$, $\psi\in$
$H^{\infty}$: We will show that if $\varphi(z)=$ $z+\psi(z)$ with
$\psi\in$ $H^{\infty}$ and $\Im(\psi(z)) > \epsilon > 0$ for all
$z\in$ $\mathbb{H}$ then
\begin{equation*}
\sigma_{e}(C_{\varphi})\supseteq\{e^{izt}:t\in [0,\infty),
z\in\mathcal{R}_{\infty}(\psi)\}\cup\{0\}
\end{equation*}
 where $\mathcal{R}_{\infty}(\psi)$ is the local essential range of $\psi$
at infinity. We use the following Theorem 13 due to Axler ([2]) to
prove the above result:

\begin{theorem}[\scshape Axler's Theorem]
Let $f\in L^{\infty}$, then there is a Blaschke product $B$ and
$\varphi\in H^{\infty}+C$ so that $f=\overline{B}\varphi$
\end{theorem}

The proof of Axler's Theorem can be found in [2]. Now we state and
prove the main result of this section:

\begin{thmc}
Let $\varphi:\mathbb{H}\rightarrow$ $\mathbb{H}$ be an analytic
self-map of $\mathbb{H}$ such that $\varphi(z)=$ $z+\psi(z)$ with
$\psi\in$ $H^{\infty}(\mathbb{H})$ and $\Im(\psi(z)) > \epsilon > 0$
for all
$z\in$ $\mathbb{H}$. Then for \\
$C_{\varphi}:H^{2}(\mathbb{H})\rightarrow$ $H^{2}(\mathbb{H})$ we
have
\begin{equation*}
\sigma_{e}(C_{\varphi})\supseteq\{e^{izt}:z\in\mathcal{R}_{\infty}(\psi),t\in
[0,\infty)\}\cup\{0\}
\end{equation*}
where $\mathcal{R}_{\infty}(\psi)$ is the local essential range of
$\psi$ at infinity.
\end{thmc}

\begin{proof}
By Proposition 5 we have if $\varphi(z)=$ $z+\psi(z)$ with $\psi\in$
$H^{\infty}(\mathbb{H})$ and $\Im(\psi(z)) > \delta > 0$ for all
$z\in$ $\mathbb{H}$ then
$$C_{\varphi}\in\Psi(L^{\infty}(\mathbb{R}),C([0,\infty])).$$ Consider
the symbol map
$$\Sigma:L^{\infty}(\mathbb{R})\rightarrow\mathcal{T}(L^{\infty}(\mathbb{R}))$$
defined by $\Sigma(\psi)=$ $T_{\psi}$. Clearly $\Sigma$ is
injective. Let $\psi_{1}=$ $\varphi_{1}\bar{\varphi_{2}}$ and
$\psi_{2}=$ $\varphi_{3}\bar{\varphi_{4}}$ with $\varphi_{j}\in$
$H^{\infty}$, $j\in$ $\{1,2,3,4\}$. Then
$T_{\overline{\varphi_{2}\varphi_{4}}}=$
$T_{\bar{\varphi_{2}}}T_{\bar{\varphi_{4}}}$. This implies that
$$T_{\psi_{1}\psi_{2}}-T_{\psi_{1}}T_{\psi_{2}}=T_{\bar{\varphi_{2}}}[T_{\bar{\varphi_{4}}},T_{\varphi_{1}}]T_{\varphi_{3}}.$$
Since $L^{\infty}(\mathbb{R})$ is spanned by such $\psi_{1}$ and
$\psi_{2}$'s we have
$$T_{\psi_{1}\psi_{2}}-T_{\psi_{1}}T_{\psi_{2}}\in
com(\mathcal{T}(L^{\infty}(\mathbb{R})))$$ for all $\psi_{1}$,
$\psi_{2}\in$ $L^{\infty}(\mathbb{R})$ (for more details see [20]
pp. 345). Since $\Sigma$ is injective, $q\circ\Sigma$ is a C*
algebra isomorphism from $L^{\infty}(\mathbb{R})$ onto
$\mathcal{T}(L^{\infty}(\mathbb{R}))/com(\mathcal{T}(L^{\infty}(\mathbb{R})))$
where $$q:\mathcal{T}(L^{\infty}(\mathbb{R}))\rightarrow
\mathcal{T}(L^{\infty}(\mathbb{R}))/com(\mathcal{T}(L^{\infty}(\mathbb{R})))$$
is the quotient map. So $M(\mathcal{T}(L^{\infty}(\mathbb{R})))=$
$M(L^{\infty}(\mathbb{R}))$. Since $$K(H^{2})\subset
com(\mathcal{T}(L^{\infty}(\mathbb{R})))$$ we also have
$$M(\mathcal{T}(L^{\infty}(\mathbb{R}))/K(H^{2}))\cong
M(L^{\infty}(\mathbb{R})).$$

Hence using Power's theorem we can identify
$$M(\Psi(L^{\infty}(\mathbb{R}),C([0,\infty]))/K(H^{2}))\cong(M_{\infty}(L^{\infty}(\mathbb{R}))\times
[0,\infty])\cup (M(L^{\infty}(\mathbb{R}))\times\{\infty\})$$ where
$$M_{\infty}(L^{\infty}(\mathbb{R}))=\{x\in
M(L^{\infty}(\mathbb{R})):x|_{C(\dot{\mathbb{R}})}=\delta_{\infty}\}$$
is the fiber of $M(L^{\infty}(\mathbb{R}))$ at infinity:

Let $(x,y)\in$ $M(L^{\infty}(\mathbb{R}))\times [0,\infty)$ so that
$x|_{C(\dot{\mathbb{R}})}=$ $\delta_{t}$ with $t\neq$ $\infty$.
Choose $a\in$ $C(\dot{\mathbb{R}})$ and $\vartheta\in$
$C([0,\infty])$ having compact supports such that
$$a(t)=\vartheta(y)=1,\quad 0\leq a\leq 1,\quad
0\leq\vartheta\leq 1,\quad a(z)< 1$$ for all $z\in$
$\mathbb{R}\backslash\{t\}$ and\quad $\vartheta(w) < 1$\quad for all
$w\in$ $[0,\infty]\backslash\{y\}$. Using the same arguments as in
equation (23) we have $$\parallel
T_{a}D_{\vartheta}\parallel_{H^{2}}<1$$ which implies that
$$(x,y)\not\in M(\Psi(\mathcal{T}(L^{\infty}(\mathbb{R})),
C([0,\infty]))/K(H^{2}(\mathbb{H}))).$$

So if $(x,y)\in
M(\Psi(L^{\infty}(\mathbb{R}),C([0,\infty]))/K(H^{2}(\mathbb{H})))$
then either $y=$ $\infty$ or $x\in$
$M_{\infty}(L^{\infty}(\mathbb{R}))$.

Let $y=$ $\infty$ and $x\in$ $M(L^{\infty}(\mathbb{R}))$. Let $a\in$
$L^{\infty}(\mathbb{R})$ and $\vartheta\in$ $C([0,\infty])$ such
that $$0\leq\vartheta\leq 1,0\leq a\leq 1\quad\textrm{and}\quad
\hat{a}(x)=\vartheta(y)=1.$$ Using the same arguments as in
equations (24) and (25) together with the fact that
\begin{equation}
\lim_{w\rightarrow\infty}\parallel K S_{w}f\parallel_{2}=0
\end{equation}
for any $K\in$ $K(H^{2}(\mathbb{H}))$ and for all $f\in$
$H^{2}(\mathbb{H})$, we have
$$\parallel D_{\vartheta}T_{a}\parallel_{e}=\parallel
D_{\vartheta}T_{a}\parallel_{H^{2}}=1.$$ This implies that
$$(x,\infty)\in M(\Psi(L^{\infty}(\mathbb{R}),
C([0,\infty]))/K(H^{2}(\mathbb{H})))$$ for all $x\in
M(L^{\infty}(\mathbb{R}))$.

Now let $x\in M_{\infty}(L^{\infty}(\mathbb{R}))$ and $y\in
[0,\infty]$. Let $a\in L^{\infty}(\mathbb{R})$ and $\vartheta\in
C([0,\infty])$ such that
$$\hat{a}(x)=\vartheta(y)=1\quad\textrm{and}\quad
0\leq\vartheta\leq 1,0\leq a\leq 1.$$ Since $a\in L^{\infty}$, by
Axler's Theorem there is a Blaschke product $B$ and $b\in
H^{\infty}+C$ so that $a=\overline{B}b$. Since
$\hat{\mid\overline{B}\mid}(x)=\mid\hat{\overline{B}}(x)\mid=1$ we
have $\mid\hat{b}(x)\mid=1$. We have $M(H^{\infty}+C)\cong
M(H^{\infty})/ \mathbb{H}$ (see Corollary 6.42 of [11]), the Poisson
kernel is also asymptotically multiplicative on $H^{\infty}+C$ (see
Lemma 6.44 of [11]) and by Carleson's Corona theorem we observe that
equation (26) is also valid for $b$ for any $\varepsilon> 0$. Since
$0\leq\mid b\mid\leq 1$ this implies that there is a $w_{0} > 0$ so
that $1-\varepsilon<\mid b(w)\mid<1$ for a.e. $w$ with $\mid
w\mid>w_{0}$. After that we use the same arguments as in equations
(27) and (39) and since $\mid\overline{B}\mid=1$ a.e., we have
$$\parallel
D_{\vartheta}T_{\overline{B}b}\parallel_{e}=\parallel
D_{\vartheta}M_{\overline{B}}M_{b}\parallel_{H^{2}}=1.$$ This
implies that $$(x,y)\in M(\Psi(L^{\infty}(\mathbb{R}),
C([0,\infty]))/K(H^{2}(\mathbb{H})))$$ for all $x\in
M_{\infty}(L^{\infty}(\mathbb{R}))$.

 If $y=$
$\infty$ then for $\varphi(z)=$ $z+\psi(z)$ with $\psi\in$
$H^{\infty}$ and
  $\Im(\psi(z)) > \delta > 0$ we have
$$(x,\infty)([C_{\varphi}])=\sum_{j=0}^{\infty}\frac{1}{j!}(\hat{\tau}(x))^{j}\hat{\vartheta_{j}}(\infty)=0$$
since $\hat{\vartheta_{j}}(\infty)= 0$ for all $j$ where $\tau$ and
$\vartheta_{j}$ are as in equation (32). If $x\in$
$M_{\infty}(L^{\infty}(\mathbb{R}))$ then as in equation (33), we
have $(x,y)([C_{\varphi}])=$ $e^{i\hat{\psi}(x)y}$. By equations
(34) and (3) we have
$$\sigma_{e}(C_{\varphi})\supseteq\{x([C_{\varphi}]):x\in
M(\Psi(L^{\infty}(\mathbb{R}),
C([0,\infty]))/K(H^{2}(\mathbb{H})))\}$$ and this implies that
$$\sigma_{e}(C_{\varphi})\supseteq\{e^{i\hat{\psi}(x)t}:x\in
M_{\infty}(L^{\infty}(\mathbb{R})),t\in [0,\infty)\}\cup\{0\}.$$ By
Proposition 10 we have
$$\{\hat{\psi}(x):x\in
M_{\infty}(L^{\infty}(\mathbb{R}))\}=\mathcal{R}_{\infty}(\psi)$$
where $\mathcal{R}_{\infty}(\psi)$ is the local essential range of
$\psi$ at infinity. Hence we have
$$\sigma_{e}(C_{\varphi})\supseteq\{e^{izt}:z\in\mathcal{R}_{\infty}(\psi),t\in
[0,\infty)\}\cup\{0\}$$
\end{proof}

In the above theorem we do not have in general equality since
$\Psi(L^{\infty}(\mathbb{R}),C([0,\infty]))$ is not commutative. And
we also have the corresponding result for the unit disc:

\begin{thmd}
Let $\varphi:\mathbb{D}\rightarrow$ $\mathbb{D}$ be an analytic
self-map of $\mathbb{D}$ such that
$$\varphi(z)=\frac{2iz+\eta(z)(1-z)}{2i+\eta(z)(1-z)}$$ where
$\eta\in$ $H^{\infty}(\mathbb{D})$ with $\Im(\eta(z)) > \epsilon> 0$
for all $z\in$ $\mathbb{D}$. Then for
$C_{\varphi}:H^{2}(\mathbb{D})\rightarrow$ $H^{2}(\mathbb{D})$ we
have
$$\sigma_{e}(C_{\varphi})\supseteq\{e^{izt}:t\in [0,\infty),
z\in\mathcal{R}_{1}(\eta)\}\cup\{0\}$$ where $\mathcal{R}_{1}(\eta)$
is the local essential range of $\eta$ at $1$.
\end{thmd}

\begin{proof}
    Repeating the same arguments as in the proof of Theorem B, we have
    $$\Phi\circ C_{\varphi}\circ\Phi^{-1}\in\Psi(L^{\infty}(\mathbb{R}),C([0,\infty]))/K(H^{2}(\mathbb{H})).$$

    Take $(x,\infty)\in M(\Psi(L^{\infty}(\mathbb{R}),
    C([0,\infty]))/K(H^{2}(\mathbb{H})))$, since
    $(x,\infty)([C_{\tilde{\varphi}}])= 0$ we have
    $$(x,\infty)([\Phi\circ C_{\varphi}\circ\Phi^{-1}])=0$$ For
    $(x,y)\in$ $M_{\infty}(L^{\infty}(\mathbb{R}))\times [0,\infty]$ we
    have $(x,y)([T_{\frac{1}{i+z}}])= 0$ and hence we have
    $$(x,y)([\Phi\circ C_{\varphi}\circ\Phi^{-1}])=(x,y)([C_{\tilde{\varphi}}])$$
    for all $x\in$ $M_{\infty}(L^{\infty}(\mathbb{R}))$. Therefore by
    equation (5) and Theorem C we have
    \begin{equation*}
 \sigma_{e}(C_{\varphi})=\sigma_{e}(\Phi\circ C_{\varphi}\circ\Phi^{-1})\supseteq\{e^{izt}:t\in [0,\infty),z\in\mathcal{R}_{1}(\eta)=\mathcal{R}_{\infty}(\eta\circ\mathfrak{C})\}\cup\{0\}
    \end{equation*}
\end{proof}

\section{Acknowledgements}
The author would like to express his deep gratitude to his post-doc
adviser Prof. Thomas L. Kriete for his help and support throughout
his visit to UVA in 2009. The author wishes to express his gratitude
to his Ph.D. adviser Prof. Ayd{\i}n Aytuna for very useful and
fruitful discussions on the subject, and to his co-adviser Prof.
Theodore W. Gamelin for his support during his visits to UCLA in
2005 and 2006. The author wishes to express his thanks to the
following people: Operator Theory group of UVA, Prof. William T.
Ross, T\"{u}rker \"{O}zsar{\i}, Mrs. Julie Riddleberger, Dr. Erdal
Karap{\i}nar, Dr. Ayla Ayalp Ross and Dale Ross. This work was
supported by grants from TUBITAK (The Scientific and Technological
Research Council of Turkey)both by a doctoral scholarship and a
post-doctoral scholarship.

% ----------------------------------------------------------------
\bibliographystyle{amsplain}

\end{document}